\magnification=1202

\input AMSSYM.def
 at9.98pt
\font\bb=msbm10 at9.98pt
\font\bfsl=cmbxsl10 at9.98pt
 at9.98pt 
\font\cmbxscript=cmbx6 at5pt
 at9.98pt
\font\eightrm=cmr8
\font\eightsc=cmcsc8
 at8.3pt
\font\labf=cmbx10 at13.1pt
 at15.74pt
\font\larm=cmr10 at13.1pt
\font\rsfs=rsfs10 at9.98pt
\font\sc=cmcsc10 at9.98pt
\font\scriptygoth=ygoth at7pt
\font\srsfs=rsfs7 at6.98pt
\font\tenpbf=cmbx10 at8.32pt
\font\tenpit=cmti10 at8.32pt
\font\tenprm=cmr10 at8.32pt
\font\ygoth=ygoth at9.98pt

\textfont4=\ygoth \scriptfont4=\scriptygoth
\def\ygot{\fam4 \ygoth}
\textfont6=\rsfs \scriptfont6=\srsfs
\def\script{\fam6 \rsfs}



\def\Ad{{\rm Ad}}

\def\an{\raise0.5pt\hbox{$\kern2pt\scriptstyle\in\kern2pt$}}
\def\Ann{\hbox{\rm Ann\kern1pt}}
\def\Anns{\hbox{$\scriptstyle\rm Ann\kern0.5pt$}}

\def\arkef{\advance\chapternumber by 1\sc\roman{\the\chapternumber}}
\def\Aut{\hbox{\rm Aut\kern1pt}}
\def\bell{\hskip0pt\lower1.6pt\hbox{\bel\char'012}\kern5pt}

\def\callige{\hbox{\calligl e\kern2pt}}
\def\cheridexi{\hskip0pt\lower2pt\hbox{\cheridexia}\kern5pt}
\def\cheridexia{{\bbding\char'21}}

\def\Coad{{\rm Coad}}
\def\coker{{\rm coker\kern1pt}}
\def\Colon{\colon\kern2pt}
\def\comp{\hbox{\lower5.8pt\hbox{\larm\char'027}}}
\def\complex{\hbox{\cmbxscript\char'103}}
\def\corang{{\rm corang\kern1pt}}
\def\cos{\hbox{\rm cos\kern1pt}}
\def\cosh{\hbox{\rm cosh\kern1pt}}
\def\dbaraux{\hbox{\= {\kern-2pt\= {}}}}
\def\dbar#1{\raise3pt\hbox{\dbaraux}\kern-7.8pt #1}

\def\dif{\hbox{$C^{\infty}$}}
\def\dim{{\rm dim\kern1pt}}
\def\double{\hbox{\kern1.5pt\bb\char'156\kern-7.6pt\char'157\kern1.5pt}}

\def\enwsh#1{{\lower2.1pt\hbox{$\buildrel{\textstyle\cup}\over
{\lower.8pt\hbox{${}_{\scriptscriptstyle#1}$}}$}}}
\def\exp{{\rm exp\kern1pt}}

\def\exten{\hbox{\callig \kern-2.5pt Ext\lower2.5pt\hbox{\kern2.5pt}}}
\def\Ham{\hbox{\rm Ham\kern1pt}}
\def\im{{\rm im\kern1pt}}
\def\k{\raise0.25pt\hbox{$\ygot k$}}

\def\ker{{\rm ker\kern1pt}}

\def\Lie{\hbox{
\callig Lie\kern2pt}}
\def\mavrodexi{\hskip0pt\lower2pt\hbox{\mavrodexia}\kern5pt}
\def\mavrodexia{{\bbding\char'15}}
\def\meriki{\hbox{\cyr\char'144\kern0.3pt}}

\def\na{\raise0.5pt\hbox{$\kern2pt\scriptstyle\ni\kern2pt$}}
\def\noan{\hbox{$\an\raise0.6pt\hbox{$\kern-6.5pt\scriptstyle
          \slash\kern3pt$}$}}
\def\odot{\;{\mathchar"220C}\;}

\def\oplus{\;{\mathchar"2208}\;}
\def\opluss#1{\hbox{$\oplus_{\hskip-0.1truecm#1\hskip0.14truecm}$}}
\def\otimes{\;{\mathchar"220A}\;}
\def\pounds{\rlap{\lower3.5pt\hbox{\kern2.9pt\hbox{\char'26}}}
           {\script L}}
\def\pr{\hbox{\kern3pt{\calligs p}\callig r\kern2pt}}
\def\qed{\hbox{\kern0.3cm\vrule height5pt width5pt depth-0.2pt}}
\def\QED{\hbox{\kern0.3cm\vrule height6pt width6pt depth-0.2pt}}
\def\R{\hbox{\bf\char'122}}
\def\rang{{\rm rang\kern1pt}}

\def\san{\raise0.5pt\hbox{$\kern0.7pt\scriptscriptstyle
         \in\kern0.7pt$}}

\def\scomp{\hskip-0.05truecm\hbox{\lower5pt\hbox{$\mathchar"2017$}}
           \hskip-0.05truecm}

\def\sem{\hbox{{\script S}\kern-2.5pt\callig em\kern2pt}}
\def\semidir{\hbox{$\;$\bb\char'156$\;$}}
\def\sin{\hbox{\rm sin\kern1pt}}
\def\sinh{\hbox{\rm sinh\kern1pt}}

\def\Sphere{\hbox{\bf\char'123}}

\def\styl{\hbox{\bbding\char'26}}
\def\stylo{\hskip0.3truecm\hbox{\lower1.5pt\hbox{\styl}}}
\def\times{\;{\mathchar"2202}\;}
\def\timess#1{\hbox{$\times_{\hskip-0.1truecm#1\hskip0.14truecm}$}}
\def\tonos{\hbox{\kern-1.3pt\lower0.7pt\hbox{$\mathchar"6013$}}}
\def\tonoskef{\hbox{$\kern-1.3pt\mathchar"6013$}}

\def\wbaraux{\hbox{\= {\kern-1.4pt\= {\kern-1.4pt\= {\kern-1.4pt\=
 {\kern-1.4pt\= {\kern-1.4pt\= {\kern-1.4pt\= {\kern-1.4pt\= {}}}}}}}}}}
\def\wbar#1{\hbox{\raise3pt\hbox{\wbaraux}\kern-30.5pt #1}}
\def\wh#1{\widehat{#1}}
\def\wt#1{\widetilde{#1}}

\def\wwbaraux{\hbox{\= {\kern-1.4pt\= {\kern-1.4pt\= {\kern-1.4pt\=
{\kern-1.4pt\= {\kern-1.4pt\= {\kern-1.4pt\= {\kern-1.4pt\=
{\kern-1.4pt\= {}}}}}}}}}}}
\def\wwbar#1{\hbox{\raise3pt\hbox{\wwbaraux}\kern-34pt #1}}


\catcode`\@=11
\def\eightpoint{\eightrm}
\def\footnote#1{\edef\@sf{\spacefactor\the\spacefactor}#1\@sf
     \insert\footins\bgroup\eightpoint
     \interlinepenalty100 \let\par=\endgraf
      \leftskip=0pt \rightskip=0pt
      \splittopskip=10pt plus 1pt minus 1pt \floatingpenalty=20000
      \smallskip\item{#1}\bgroup\strut\aftergroup\@foot\let\neft}
\skip\footins=12pt plus 2pt minus 4pt
\dimen\footins=30pc

\def\line{\hbox to\hsize}

\def\title#1{\line{\hss}\line{\hss#1\hss}%
\line{\hss}\hskip-0.75truecm}

\def\author#1{{\tenprm #1:}}
\def\ekdoths#1{{\tenprm #1}}
\def\periodiko#1{{\tenpit #1\tenprm ,}}
\def\selides#1{{\tenprm #1}}
\def\titlosa#1{{\tenprm #1,}}
\def\titlosb#1{{\tenpit #1\tenprm ,}}
\def\volume#1{{\tenprm Vol. \tenpbf #1\tenprm :}}

%
%
%
\def\teleia{\hbox{.}}
\newif\ifPhysRev
\def\Textindent#1{\noindent\llap{#1\enspace}\ignorespaces}
\def\nonfrenchspacing{\sfcode`\.=3001 \sfcode`\!=3000 \sfcode`\?=3000
        \sfcode`\:=2000 \sfcode`\;=1500 \sfcode`\,=1251 }
\nonfrenchspacing
\newdimen\d@twidth
 {\setbox0=\hbox{s.} \global\d@twidth=\wd0 \setbox0=\hbox{s}
        \global\advance\d@twidth by -\wd0 }
\def\removehglue{\loop \unskip \ifdim\lastskip >\z@ \repeat }
\def\roll@ver#1{\removehglue \nobreak \count255 =\spacefactor \dimen@=\z@
        \ifnum\count255 =3001 \dimen@=\d@twidth \fi
        \ifnum\count255 =1251 \dimen@=\d@twidth \fi
    \iftwelv@ \kern-\dimen@ \else \kern-0.83\dimen@ \fi
   #1\spacefactor=\count255 }
\def\step@ver#1{\relax \ifmmode #1\else \ifhmode
        \roll@ver{${}#1$}\else {\setbox0=\hbox{${}#1$}}\fi\fi }
\def\attach#1{\step@ver{\strut^{\mkern 2mu #1} }}

\normalbaselineskip = 20pt plus 0.2pt minus 0.1pt
\normallineskip = 1.5pt plus 0.1pt minus 0.1pt
\normallineskiplimit = 1.5pt
\newskip\normaldisplayskip
\normaldisplayskip = 20pt plus 5pt minus 10pt
\newskip\normaldispshortskip
\normaldispshortskip = 6pt plus 5pt
\newskip\normalparskip
\normalparskip = 6pt plus 2pt minus 1pt
\newskip\skipregister
\skipregister = 5pt plus 2pt minus 1.5pt
\newif\ifsingl@    \newif\ifdoubl@
\newif\iftwelv@    \twelv@true
\def\singlespace{\singl@true\doubl@false\spaces@t}
\def\doublespace{\singl@false\doubl@true\spaces@t}
\def\normalspace{\singl@false\doubl@false\spaces@t}
\def\Tenpoint{\tenpoint\twelv@false\spaces@t}
\def\Twelvepoint{\twelvepoint\twelv@true\spaces@t}
\def\spaces@t{\relax
      \iftwelv@ \ifsingl@\subspaces@t3:4;\else\subspaces@t1:1;\fi
       \else \ifsingl@\subspaces@t3:5;\else\subspaces@t4:5;\fi \fi
      \ifdoubl@ \multiply\baselineskip by 5
         \divide\baselineskip by 4 \fi }
\def\subspaces@t#1:#2;{
      \baselineskip = \normalbaselineskip
      \multiply\baselineskip by #1 \divide\baselineskip by #2
      \lineskip = \normallineskip
      \multiply\lineskip by #1 \divide\lineskip by #2
      \lineskiplimit = \normallineskiplimit
      \multiply\lineskiplimit by #1 \divide\lineskiplimit by #2
      \parskip = \normalparskip
      \multiply\parskip by #1 \divide\parskip by #2
      \abovedisplayskip = \normaldisplayskip
      \multiply\abovedisplayskip by #1 \divide\abovedisplayskip by #2
      \belowdisplayskip = \abovedisplayskip
      \abovedisplayshortskip = \normaldispshortskip
      \multiply\abovedisplayshortskip by #1
        \divide\abovedisplayshortskip by #2
      \belowdisplayshortskip = \abovedisplayshortskip
      \advance\belowdisplayshortskip by \belowdisplayskip
      \divide\belowdisplayshortskip by 2
      \smallskipamount = \skipregister
      \multiply\smallskipamount by #1 \divide\smallskipamount by #2
      \medskipamount = \smallskipamount \multiply\medskipamount by 2
      \bigskipamount = \smallskipamount \multiply\bigskipamount by 4 }
\def\normalbaselines{ \baselineskip=\normalbaselineskip
   \lineskip=\normallineskip \lineskiplimit=\normallineskip
   \iftwelv@\else \multiply\baselineskip by 4 \divide\baselineskip by 5
     \multiply\lineskiplimit by 4 \divide\lineskiplimit by 5
     \multiply\lineskip by 4 \divide\lineskip by 5 \fi }


\def\abstract#1{\parshape=1 0.7cm \dimen10
                {\tenpbf Abstract. \tenprm #1}}

\newcount\appendixnumber     \appendixnumber=0
\newcount\chapternumber      \chapternumber=0
\newcount\equanumber         \equanumber=0
\newcount\mathnumber         \mathnumber=0
\newcount\appequanumber      \appequanumber=0
\newcount\appmathnumber      \appmathnumber=0

\let\variableone=\relax
\let\variabletwo=\relax
\let\chapterlabel=\relax
\let\sectionlabel=\relax
\let\mathlabel=\relax
\newtoks\chapterstyle        \chapterstyle={\Number}
\newtoks\sectionstyle        \sectionstyle={\chapterlabel\Number}
\newskip\chapterskip         \chapterskip=\bigskipamount
\newskip\sectionskip         \sectionskip=\medskipamount
\newskip\headskip            \headskip=8pt plus 3pt minus 3pt
\newdimen\chapterminspace    \chapterminspace=15pc
\newdimen\sectionminspace    \sectionminspace=10pc
\newdimen\sectionspace       \sectionspace=20pc
\newdimen\referenceminspace  \referenceminspace=25pc

\def\chapterreset{\global\advance\chapternumber by 1
   \ifnum\equanumber<0 \else\global\equanumber=0\fi
   \mathnumber=0
   \makechapterlabel}
\def\makechapterlabel{\let\sectionlabel=\relax\let\mathlabel=\relax
 \xdef\chapterlabel{\the\chapterstyle{\the\chapternumber\teleia\kern3pt}}}

\def\rightheadline{\sc\hfil\variableone\eightsc\hfil\folio}
\def\leftheadline{\eightsc\folio\hfil{\sc\variabletwo}\hfil}
\def\heads{\footline={\hfil}\headline={\ifodd\pageno
               \rightheadline\else\leftheadline\fi}}

\def\headseis{\partreset\headline={\ifodd\pageno{
                         \hfil\sc partie {\eightsc\the\partnumber}
                         -introduction\hfil\eightsc\folio}\else
                        {\eightsc\folio\hfil\sc partie
                         {\eightsc\the\partnumber}-introduction\hfil}\fi}
                        \footline={\hfil}}

\def\alphabetic#1{\count255='140 \advance\count255 by #1\char\count255}
\def\Alphabetic#1{\count255='100 \advance\count255 by #1\char\count255}
\def\Roman#1{\uppercase\expandafter{\romannumeral #1}}
\def\roman#1{\romannumeral #1}
\def\Number#1{\number #1}
\def\BLANC#1{}

\def\titlestyle#1{\par\begingroup \interlinepenalty=9999
     \leftskip=0.02\hsize plus 0.23\hsize minus 0.02\hsize
     \rightskip=\leftskip \parfillskip=0pt
     \hyphenpenalty=9000 \exhyphenpenalty=9000
     \tolerance=9999 \pretolerance=9000
     \spaceskip=0.333em \xspaceskip=0.5em
     \iftwelv@\bf\else\bf\fi
   \noindent #1\par\endgroup }

\def\spacecheck#1{\dimen@=\pagegoal\advance\dimen@ by -\pagetotal
   \ifdim\dimen@<#1 \ifdim\dimen@>0pt \vfil\break \fi\fi}
\def\TableOfContentEntry#1#2#3{\relax}

\def\chapter#1{\par\vskip0.7cm
   \chapterreset \titlestyle{\chapterlabel\ #1}
   \nobreak\vskip\headskip
   \wlog{\string\chapter\space \chapterlabel} }

\def\appendixreset{\global\advance\appendixnumber by 1
                   \appmathnumber=0\appequanumber=0}
\def\appendix#1{\par \penalty-300\vskip\chapterskip
   \spacecheck\chapterminspace
   \appendixreset \title{\bf Appendix \Alphabetic{\the\appendixnumber}}
   \nobreak\vskip-\chapterskip\penalty 30000
   \vskip-\chapterskip
   \par{\titlestyle{#1}}
   \vskip\chapterskip
   \wlog{\string\appendix\space \chapterlabel} }

%
%
\def\eqname#1{\relax \ifnum\equanumber<0
     \xdef#1{{\noexpand\rm(\number-\equanumber)}}%
       \global\advance\equanumber by -1
    \else \global\advance\equanumber by 1
      \xdef#1{{\noexpand(\rm{\the\chapternumber}\teleia
                            \rm{\number\equanumber})}} \fi #1}

\def\eqn{\eqno\eqname}

\def\math#1#2{\vskip0.1cm
   \global\advance\mathnumber by 1
   \xdef\mathlabel{\the\chapternumber\teleia\the\mathnumber}
   \wlog{\string\math\space \mathlabel}
   {\bf\enspace\mathlabel\hskip0.2cm #1}
   \xdef#2{{\mathlabel}}}

\def\appeqname#1{\relax \ifnum\appequanumber<0
     \xdef#1{{\noexpand\rm(\number-\appequanumber)}}%
       \global\advance\appequanumber by -1
    \else \global\advance\appequanumber by 1
      \xdef#1{{\noexpand(\hbox{\Alphabetic{\the\appendixnumber}}\teleia
                            {\number\appequanumber})}} \fi #1}

\def\mathapp#1#2{\vskip0.1cm
   \global\advance\appmathnumber by 1
   \xdef\appmathlabel{{\Alphabetic{\the\appendixnumber}}\teleia
   \the\appmathnumber}
   \wlog{\string\mathapp\space \appmathlabel}
   {\bf\enspace\appmathlabel\hskip0.2cm #1}
   \xdef#2{{\appmathlabel}}}


%
%
%
\newtoks\referencestyle      \referencestyle={\tenpbf\Number}
\newcount\referencecount     \referencecount=0
\newcount\lastrefsbegincount \lastrefsbegincount=0
\newif\ifreferenceopen       \newwrite\referencewrite
\newif\ifrw@trailer
\newdimen\refindent     \refindent=13pt
\def\NPrefmark#1{\attach{\scriptscriptstyle [ #1 ] }}
\let\PRrefmark=\attach
\def\refmark#1{\relax\ifPhysRev\PRrefmark{#1}\else\NPrefmark{#1}\fi}
\def\refend@{\refmark{\number\referencecount}}
\def\refend{\refend@{}\space }
\def\refsend{\refmark{\count255=\referencecount
   \advance\count255 by-\lastrefsbegincount
   \ifcase\count255 \number\referencecount
   \or \number\lastrefsbegincount,\number\referencecount
   \else \number\lastrefsbegincount-\number\referencecount \fi}\space }
\def\refitem#1{\par\hangafter=0 \hangindent=\refindent	\Textindent{#1}}
\def\Ref{\rw@trailertrue\REF}
\def\REF#1{\r@fstart{#1}%
   \rw@begin{\tenprm [\tenpbf\Number{\the\referencecount}\tenprm ]}\rw@end}
\def\r@fstart#1{\chardef\rw@write=\referencewrite \let\rw@ending=\refend@
   \ifreferenceopen \else \global\referenceopentrue
   \immediate\openout\referencewrite=referenc.txa
   \toks0={\catcode`\^^M=10}\immediate\write\rw@write{\the\toks0} \fi
   \global\advance\referencecount by 1 
   \xdef#1{[{\the\referencestyle{\the\referencecount}}]}}
 {\catcode`\^^M=\active %
 \gdef\rw@begin#1{\immediate\write\rw@write{\noexpand\refitem{#1}}%
   \begingroup \catcode`\^^M=\active \let^^M=\relax}%
 \gdef\rw@end#1{\rw@@end #1^^M\rw@terminate \endgroup%
   \ifrw@trailer\rw@ending\global\rw@trailerfalse\fi }%
 \gdef\rw@@end#1^^M{\toks0={#1}\immediate\write\rw@write{\the\toks0}%
   \futurelet\n@xt\rw@test}%
 \gdef\rw@test{\ifx\n@xt\rw@terminate \let\n@xt=\relax%
       \else \let\n@xt=\rw@@end \fi \n@xt}%
}
\let\rw@ending=\relax
\let\rw@terminate=\relax

\def\Closeout#1{\toks0={\catcode`\^^M=5}\immediate\write#1{\the\toks0}%
   \immediate\closeout#1}

\input pictex
\topskip 1truecm
\voffset=2.5truecm
\hsize 15truecm
\vsize 21truecm
\hoffset=0.25cm
\def\undertext#1{$\underline{\hbox{#1}}$}
\baselineskip=14pt plus .2pt
\topglue 3truecm
\pageno=1
\def\variableone{semidirect products and the pukanszky condition}
\def\variabletwo{p. baguis}
\dimen10=\hsize \advance\dimen10 by -1.4cm
\footline={\hss\eightsc\folio\hss}

\Ref\dbkp{
\author{Duval, C., Burdet, G., K\"unzle, H. P., Perrin, M.}
\titlosa{Bargmann structures and Newton-Cartan theory} 
\periodiko{Phys. Rev. D}
\volume{31(8)}
\selides{1841--1853 (1985)}}

\Ref\de{
\author{Duval, C., Elhadad, J.}
\titlosa{Geometric quantization and localization of relativistic spin
systems} 
\periodiko{Contemp. Math.}
\volume{132}
\selides{317--330 (1992)}}

\Ref\det{
\author{Duval, C., Elhadad, J., Tuynman, G.M.}
\titlosa{Pukanszky's condition and symplectic induction} 
\periodiko{J. Diff. Geom.}
\volume{36}
\selides{331--348 (1992)}}

\Ref\guist{
\author{Guillemin, V., Sternberg, S.}
\titlosb{Symplectic techniques in physics}
\ekdoths{Cambridge University Press, Cambridge (1991)}}

\Ref\kazkost{
\author{Kazhdan, D., Kostant, B., Sternberg, S.}
\titlosa{Hamiltonian Group Actions and Dynamical Systems of Calogero
Type}
\periodiko{Comm. Pure Appl. Math.}
\volume{31}
\selides{481--508 (1978)}}

\Ref\kob{
\author{Kobayashi, S., Nomizu, K.}
\titlosb{Foundations of Differential Geometry}
\ekdoths{Vol. I, II, Interscience Publishers, 
John Wiley \& Sons, (1963, 1969)}}

\Ref\kost{
\author{Kostant, B.}
\titlosa{On certain unitary representations which arise from a
quantization theory}
\periodiko{in ``Group representations in mathematics and physics"}
\selides{(V. Bargman, ed.), pp. 237--253, 
Proceedings Battelle Seattle 1969 Rencontres; LNP {\tenpbf 6},}
\selides{Springer-Verlag, Berlin, 1970}}

\Ref\land{
\author{Landsman, N. P.}
\titlosa{Rieffel induction as generalized quantum Marsden-Weinstein
reduction}
\periodiko{J. Geom. Phys.}
\volume{15}
\selides{285--319 (1995)}}

\Ref\maraw{
\author{Marsden, J. E., Ratiu, T., Weinstein, A.}
\titlosa{Semidirect products and reduction in mechanics}
\periodiko{Trans. A.M.S}
\volume{281}
\selides{147--177 (1984)}}

\Ref\pukanszky{
\author{Pukanszky, L.}
\titlosa{On the theory of exponential groups}
\periodiko{Trans. A.M.S}
\volume{126}
\selides{487--507 (1967)}}

\Ref\ra{
\author{Rawnsley, J. H.}
\titlosa{Representations of a semidirect product by quantization} 
\periodiko{Math. Proc. Camb. Phil. Soc.}
\volume{78}
\selides{345--350 (1975)}}

\Ref\sn{
\author{\'Sniatycki, J.}
\titlosb{Geometric quantization and quantum mechanics}
\ekdoths{Appl. Math. Sci.,}
\volume{50}
\ekdoths{Springer-Verlag (1980)}}

\Ref\sou{
\author{Souriau, J.-M.}
\titlosb{Structure des syst\`emes dynamiques} 
\ekdoths{Dunod, Paris (1969)}}

\Ref\stern{
\author{Sternberg, S.}
\titlosa{Minimal coupling and the symplectic mechanics of a classical
particle in the presence of a Yang-Mills field}
\periodiko{Proc. Nat. Acad. Sci.}
\volume{74}
\selides{5253--5254 (1977)}}

\Ref\wein{
\author{Weinstein, A.}
\titlosa{A universal phase space for particles in Yang-Mills fields}
\periodiko{Lett. Math. Phys.}
\volume{2}
\selides{417--420 (1978)}}

\Ref\zak{
\author{Zakrzewski, S.}
\titlosa{Induced representations and induced Hamiltonian actions}
\periodiko{J. Geom. Phys.}
\volume{3}
\selides{211--219 (1986)}}

\Ref\zie{
\author{Ziegler, F.}
\titlosb{M\'ethode des orbites et repr\'esentations quantiques}
\ekdoths{Th\`ese, Univesrit\'e de Provence, D\'ecembre 1996}}


\advance\baselineskip by -2.5pt

\titlestyle{\labf Semidirect products and the Pukanszky condition}

\vskip0.5cm

\centerline{\bf P. Baguis}

\vskip0.3cm

\centerline{Centre de Physique Th\'eorique-CNRS}
\centerline{Luminy, Case 907-F-13288}
\centerline{Marseille Cedex 9, France}

\centerline{e-mail: {\tt baguis@cpt.univ-mrs.fr}}

\vskip0.5cm
\vskip0.5cm

\abstract{We study the general geometrical structure of the coadjoint
orbits of a semidirect product formed by a Lie group and a
representation of this group on a vector space.  
The use of symplectic induction methods gives
new insight into the structure of these orbits. In fact,
each coadjoint orbit of such a group is
obtained by symplectic induction on some coadjoint orbit of a
``smaller" Lie group. We study also a
special class of polarizations related to a semidirect product
and the validity of Pukanszky's condition for these polarizations.
Some examples of physical interest are discussed using the previous
methods.}

\advance\baselineskip by 2.5pt
 
\vskip1cm
 
{\tenpit Key-words:} {\tenprm Lie groups, coadjoint orbits,
Pukanszky's condition, symplectic induction}
 
\vskip0.1cm
 
{\tenpit 1991 MSC:} {\tenprm 53C15}

\vskip0.1cm

{\tenpit 1996 PACS:} {\tenprm 02.20.Qs}

\vskip0.5cm

\chapter{Introduction}

Polarized coadjoint orbits of a Lie group $G$, are good
candidates for geometrically quantized phase spaces. They also play a
central r\^ole in representation theory, more specifically in the
context of Kirillov's ``orbit method". In the case of exponential
groups, Pukanszky showed {\pukanszky} that the orbit method leads to irreducible
unitary representations of $G$, if and only if the polarization
satisfies a certain condition, known as Pukanszky's condition (see
Lemma 6.2 (1) with $\don=\frak e$). This method has been adapted to
the case of complex polarizations, especially for solvable Lie groups:
Auslander and Kostant {\kost}, showed that Pukanszky's condition was
needed  in order to guarantee the irreducibility of the representations
obtained via holomorphic induction from the real polarizing subgroup
$D\subset G$ (see below).

In its initial formulation, Pukanszky's condition means that the coadjoint
orbit in question contains an affine plane, constructed out by the
polarization. Only recently {\de} {\det},  has it been realized that
validity of Pukanszky's condition is equivalent to the fact that the
corresponding coadjoint orbit is symplectomorphic to a modified
cotangent bundle, obtained by symplectic induction from a point. The 
physical consequences of this symplectomorphism have been studied in the
previous references for the coadjoint orbits of the Poincar\'e
group, which is a semidirect product. 

\heads

Our aim is to give, on the one hand, a detailed analysis of the
geometrical structure of the semidirect product coadjoint orbits,
for the case where this product is formed by a Lie group $K$ and
a representation $\rho\Colon K\rightarrow GL(V)$ on a vector
space $V$, {\ra}. On the other hand, we consider this very interesting
geometrical structure in the framework of Pukanszky's condition.
Summarizing the results of this article, we mention the following
three points:
\item{$\bullet$} the coadjoint orbits of a semidirect product present
several analogies with the cotangent bundles and under certain
conditions they are in reality cotangent bundles of $K$-orbits
in the dual $V^{\ast}$ of the vector space $V$;
\item{$\bullet$} the validity of Pukanszky's condition for a
special class of polarizations of the semidirect product $G=
K\timess{\rho}V$ is equivalent to the validity of the same condition
for ``smaller" polarizations associated to the homogeneous part
$K$;
\item{$\bullet$} the coadjoint orbits of the semidirect product
$G=K\timess{\rho}V$ are obtained by symplectic induction
on coadjoint orbits of appropriate subgroups of the homogenous part
$K$.

\noindent In what concerns the third point, a variant of the
symplectic induction we are using here, gave recently
the same result, {\land}; see also {\zie} for a more general treatment
in the context of symplectic Mackey's theory.

We discuss finally three examples of semidirect product
and we apply the previous formalism in the geometrical study of their
coadjoint orbits. The semidirect products in question are the special
Euclidean, the Galilei and the Bargmann group,
three Lie groups whose coadjoint orbits are respectively related to
geometrical optics, polarization of light and to the dynamics of
non-relativistic particles.

\vskip0.3cm

{\bf Acknowledgments.} I would like to thank Professor C. Duval for
careful and critical reading of the manuscript as well as for useful
and stimulating discussions. In addition, I thank F. Ziegler for
having first pointed out that semidirect product coadjoint orbits are
indeed induced symplectic manifolds. 
Let me also thank Professor J. H. Rawnsley for his kind
interest in this work.

\chapter{The semidirect product}

In this section, we fix the notation concerning the semidirect
product following {\ra}.

Consider a Lie group $K$ with Lie algebra $\k$; let 
$(\ek,f)\mapsto\ek\cdot f$ be the coadjoint representation of $K$ on 
$\k^{\ast}$, the dual of the Lie algebra $\k$, and $(A,f)\mapsto A\cdot f$
its derivative, $\ek\an K,A\an\k,f\an\k^{\ast}$. If 
$\rho\Colon K\rightarrow GL(V)$ is a representation
of $K$ on the vector space $V$, then we note $\rho(\ek)v=\ek\cdot v$,
$\ek\an K,v\an V$.
We note accordingly by $(\ek,p)\mapsto\ek\cdot p$ the 
contragredient representation of $\rho$, $p\an V^{\ast}$ 
and by $(A,v)\mapsto A\cdot v$ and $(A,p)\mapsto A\cdot p$ 
the corresponding derivative representations of $\k$ on $V$ and 
$V^{\ast}$.

We form now the semidirect product  $G=K
\timess{\rho}{\mit V}$. As a set $G=K\times V$ and the
group operation in $G$ is given by $$(\ek,v)\cdot
(\el,u)=(\ek\el,\ek\cdot u+v),\;\forall (\ek,v),(\el,u)\an G.
\eqn\groupoper$$ 
When the representation $\rho$ is understood, we write simply
$G=K\semidir V$.

The Lie algebra of this group is ${\frak g}=\k\oplus V$ (as a vector
space)
and the Lie algebra structure is given by the bracket 
$$[(A,a),(B,b)]=([A,B],A\cdot b-B\cdot a),
\;\forall (A,a),(B,b)\an{\frak g}.\eqn\commut$$
We will note $\frak g=\k\opluss{\rho}V$.

By identifying the dual $\frak g^{\ast}$ of $\frak g$ with 
$\k^{\ast}\oplus V^{\ast}$, we can express the duality between $\frak
g$ and $\frak g^{\ast}$ as
$$\em(\xi)=f(A)+p(a),\;\forall\em=(f,p)\an\frak g^{\ast},
\xi=(A,a)\an\frak g\eqn\duality$$
and the adjoint and coadjoint representations 
of $G$ on $\frak g$ and 
$\frak g^{\ast}$ respectively, by the following relations:
$$\Ad(\ek,v)(A,a)=\big(\Ad(\ek)A,\ek\cdot a-\rho\tonos\big(\Ad
(\ek)A\big)v\big),\;\forall (\ek,v)\an G,
(A,a)\an\frak g,\eqn\adjoint$$ 
$$\Coad(\ek,v)(f,p)=(\ek\cdot f+\ek\cdot p\odot v,\ek\cdot p), 
\;\forall (\ek,v)\an G,(f,p)
\an \frak g^{\ast}\eqn\coadjoint$$
where $p\odot v$ is the element of $\k^{\ast}$ defined by 
$$(p\odot v)(A)=p(A\cdot v)=-(A\cdot p)(v),\;\forall A\an 
\k,p\an V^{\ast},v\an V.\eqn\oginomeno$$
For $p\an V^{\ast}$ we denote by $K_{p}$ the isotropy subgroup of 
$p$ formed by those $\ek\an K$ such that
$\ek\cdot p=p$. It is clear that the Lie algebra of 
$K_{p}$ is given by the vector space $\k_{p}=\{A\an\k
\;|\;A\cdot p=0\}$. Then, if we define the linear map
$\tau_{p}\Colon\k\rightarrow {\mit V}^{\ast}$
by $$\tau_{p}(A)=-A\cdot p, \;\forall
A\an \k,\eqn\apeiktau$$ we have the equality $\k_{p}=\ker\tau_{p}$.

We express now the element $p\odot v\an\k^{\ast}$ in terms of 
the map $\tau_{p}$. The dual $\tau_{p}^{\ast}\Colon
V\rightarrow \k^{\ast}$ of $\tau_{p}$ is
given by the relation $\tau_{p}^{\ast}(v)(A)=\tau_{p}(A)(v)=
-(A\cdot p)(v)$, and so 
$\tau^{\ast}_{p}(v)=p\odot v,\;\forall p\an V^{\ast},\forall
v\an V$.

Let now $\k_{p}^{\circ}$ be the annihilator of $\k_{p}$;
then, if $i_{p}^{\ast}\Colon\k^{\ast}\rightarrow 
\k_{p}^{\ast}$ is the projection, $i_{p}\Colon\k_{p}\hookrightarrow\k$,
we have $\k_{p}^{\circ}=\ker
i_{p}^{\ast}$. The following is a useful lemma from {\ra}, giving
a characterization of the annihilator $\k_{p}^{\circ}$ in terms of the
linear map $\tau_{p}$.

\math{Lemma.}{\rangeoftau}{\sl $\k_{p}^{\circ}=\im\tau_{p}^{\ast}.$}

\undertext{\it Proof}. We give here a different and simple proof using
Lagrange multipliers. Indeed, we have
$\k_{p}=\tau_{p}^{-1}(0)\subset\k$ and the element $A\an\k_{p}$ is a
critical point of the map $i_{p}^{\ast}f\Colon\k_{p}\rightarrow\R$, for
$f\an\dif(\k)$, if and only if there exists an element $v\an
V^{\ast\ast}\cong V$ such that $A$ is a critical point of
$f-v\comp\tau_{p}$. Choosing $f\an\k^{\ast}$ (and therefore linear), we
find that $f\an\k_{p}^{\circ}$ if and only if $\exists v\an V$:
$f=v\comp\tau_{p}=p\odot v$. \QED

\chapter{Orbits-Isotropy subgroups}

We recall now the structure of the coadjoint orbits of a semidirect
product studied by Rawnsley {\ra}, and we exploit in more detail the
structure of the isotropy subgroups with respect to the coadjoint
representation for the semidirect product.
According to {\ra}, the coadjoint orbits of a semidirect product are
classified by the coadjoint orbits of ``little" groups, which are 
isotropy subgroups of its homogeneous part (see also {\guist}). In fact, 
fibre bundles having these coadjoint orbits as fibres, completely
characterize the coadjoint orbits of the semidirect product. As we
shall see later on (section 10), the little-group coadjoint orbits
play an even deeper r\^ole for the geometrical structure of the
corresponding semidirect product coadjoint orbit. 

Let now ${Z}={\script O}_{p}^{K},\;p\an V^{\ast}$ be an orbit
of $K$ in $V^{\ast}$ with respect to the representation
$\rho^{\ast}$. A bundle of little-group orbits over $Z$ is
a fibre bundle $\pi\Colon Y\rightarrow {Z}$ such that each fibre
$Y_{p}=\pi^{-1}(p)$ be a coadjoint orbit of the isotropy subgroup
$K_{p}$.

We construct the bundle of little-group orbits as
follows. Consider elements $p\an V^{\ast}$, $\phi\an\k_{p}^{\ast}$
and let $Z$, $Y_{p}$ be the corresponding orbits under the
actions of $K$ and $K_{p}$ respectively, ${Z}={\script
O}_{p}^{K},\; {Y}_{p}={\script O}_{\phi}^{K_{p}}$. There
is a left action of 
the isotropy subgroup $K_{p}\subset K$ on the product 
$K\times {Y}_{p}$ given by: $$h\cdot
(\ek,\phi)=(\ek h^{-1},h\cdot\phi).\eqn\leftaction$$
We define the bundle of little-group orbits $Y$ as the quotient 
${Y}=\bigl(K\times {Y}_{p}\bigr)\big/K_{p}$, i.e., $Y$ is
the fibre bundle associated to the principal bundle
$K\rightarrow {Z}$ with respect to the coadjoint action of
$K_{p}$ on ${Y}_{p}$. The group $K$ acts on $Y$ as
follows. If $\psi\an {Y}_{p}$, then for $\ek\an K$ we define the
point $\ek\cdot\psi\an {Y}_{\ek\cdot p}$ as 
$$(\ek\cdot\psi)(A)=\psi(\Ad(\ek^{-1})A),\;\forall A\an
\k_{\ek\cdot p}.\eqn\actionlittle$$
Consequently, the following choice to represent the points of $Y$ 
($K_{p}$-orbits in $K\times {Y}_{p}$) is appropriate:
$K_{p}\cdot (\ek,\phi)=\ek\cdot\phi\an {Y}_{\ek\cdot p}$.
We can now define the projection $\pi\Colon{Y}\rightarrow {Z}$
by $\pi\big(K_{p}\cdot
(\ek,\phi)\big)=\ek\cdot p\an {Z}$. 


It is easy to verify that this construction is independent of the
point $p\an {Z}$. Then, the  following proposition {\ra} clarifies 
the role of the bundles of little-group orbits; see also section 10 below.

\math{Proposition.}{\bundleoforbits}{\sl There is a bijection between
the set of bundles of little-group orbits and the set of coadjoint
orbits of $G$ on $\frak g^{\ast}$.}

\vskip0.3cm

Consider now the coadjoint orbit ${\script O}_{\en}^{G}$ of
$\en=(f,p)$ in $\frak g^{\ast}$ and the corresponding fibre bundle of
little-group orbits $Y$. This orbit is fibred over the $K$-orbit
$Z$ of the point $p\an V^{\ast}$ and the fibre is of the form
$({\script O}^{G}_{\en})_{q}=K_{q}\cdot h+q\odot V,\;q=\ek\cdot p\an 
{Z},\;h=\ek\cdot f$ for $\ek\an K$ (relation
$\coadjoint$). Thus if $y\an({\script O}^{G}_{\en})_{q}$,
then there exists an element $(\el,v)\an K_{q}\timess{\rho}V$
such that $y=\el\cdot h+q\odot v$ and so
$i_{q}^{\ast}(y)=i_{q}^{\ast}\big(\el\cdot h+q\odot
v\big)=i_{q}^{\ast}\big(\el\cdot h\big)=\el\cdot
i_{q}^{\ast}h\an {Y}_{q}$. It is clear that the projection
$i_{q}^{\ast}\Colon\k^{\ast}\rightarrow\k_{q}^{\ast}$ defines in fact
a projection $i_{q}^{\ast}\Colon({\script O}^{G}_{\en})_{q}
\rightarrow {Y}_{q}$
between the fibres of ${\script O}_{\en}^{G}$ and $Y$.
Furthermore, the fibre $(i_{q}^{\ast})^{-1}(\varphi)$ for $\varphi\an 
{Y}_{q}$, is the orbit of the point $q\an V^{\ast}$ under the action
of the linear subgroup $V\subset G$. We have thus proved the lemma,
{\ra}:

\math{Lemma.}{\coadorbitstruc}{\sl The coadjoint orbit ${\script O}_{\en}
^{G}$ of the element
$\en=(f,p)\an\frak g^{\ast}=\k^{\ast}\oplus V^{\ast}$ is a fibre
bundle over the bundle of little-group orbits whose typical fibre
is the orbit of $p\an V^{\ast}$ under the action of the subgroup $V\subset
{G}$.}

\vskip0.3cm

We summarize the previous results in the commutative diagram

\vskip0.5cm

\centerline{\beginpicture
\setcoordinatesystem units <0.80000cm,0.80000cm>
\setplotsymbol ({\fiverm .})
\setlinear

\arrow <1mm> [0.1,1] from 10.9 19.050 to 10.9 15.716
\arrow <1mm> [0.1,1] from 11.40 19.2 to 13.652 17.621
\arrow <1mm> [0.1,1] from 13.74 17.25 to 11.26 15.558

\put{${\script O}^{G}_{\en}$} [lB] at 10.636 19.209
\put{$Y={\script O}^{G}_{\en}/V$} [lB] at 13.811 17.304
\put{$Z$} [lB] at 10.65 15.240
\put{$\pi$} [lB] at 12.6 16.1
\put{$\Pi$} [lB] at 10.4 17.4
\endpicture
}

\vskip-0.7cm

\noindent where $\Pi\Colon{\script O}^{G}_{\en}\rightarrow{Z}$
is the projection.

We study finally the isotropy subgroup ${G}_{\en}$ of the point
$\en=(f,p)\an\frak g^{\ast}$ with respect to the coadjoint action.
Let $\phi=i_{p}^{\ast}f$ and $(K_{p})_{\phi}$ be the isotropy
subgroup of $\phi\an\k_{p}^{\ast}$ with respect to the coadjoint
action of $K_{p}$ on $\k_{p}^{\ast}$. If $(\ek,v)\an G_{\en}$, 
then $\ek\cdot p=p$ and $\ek\cdot f+p\odot v=f$
which means that $\ek\cdot\phi=\phi\Rightarrow \ek\an
(K_{p})_{\phi}$. We have thus an epimorphism $j\Colon
G_{\en}\rightarrow(K_{p})_{\phi}$ given by $j(\ek,v)=\ek$. The
kernel of $j$ is calculated easily: $\ker j=\{(\ek,v)\an
G_{\en}\,|\,j(\ek,v)=e\}=\{(e,v)\an G_{\en}\}$. But
the element $(e,v)$ belonging to $G_{\en}$ is such that
$e\cdot f+p\odot v=f\Rightarrow p\odot v=0$, thus
$v\an\ker\tau_{p}^{\ast}$. On the other hand, $\ker\tau_{p}^{\ast}$ is
a vector subgroup of $G_{\en}$ as the inclusion map
$i\Colon\ker\tau_{p}^{\ast}\rightarrow G_{\en}$ given by
$i(v)=(e,v)$ indicates.

We conclude that $\ker j=\ker\tau_{p}^{\ast}$ and we have the
following exact sequence $$0\longrightarrow
\ker\tau_{p}^{\ast}\buildrel i \over\longrightarrow 
{G}_{\en}\buildrel j \over\longrightarrow 
({K}_{p})_{\phi}\longrightarrow e\eqn\exactseq$$
which gives us all the possible information about the structure of the
isotropy subgroup $G_{\en}$.

We note that $G_{\en}$ is not in general equal to the semidirect
product of $(K_{p})_{\phi}$ and $\ker\tau_{p}^{\ast}$ because
$(K_{p})_{\phi}$ is not in general a subgroup of 
$G_{\en}$. In fact, if it were, we would have $(\ek,0)\an
G_{\en}$ for each $\ek\an(K_{p})_{\phi}$. But in such a case we
find $\ek\cdot f+p\odot 0=f\Rightarrow\ek\an K_{f}$, so 
$(K_{p})_{\phi}\subset K_{f}$; clearly this condition is not in
general satisfied.

Conversely now, the inclusion $(K_{p})_{\phi}\subset K_{f}$
induces a group monomorphism $m\Colon(K_{p})_{\phi}\rightarrow
G_{\en}$ given by $m(\ek)=(\ek,0)$ for in that case we have
$\ek\cdot f+p\odot 0=f$ for each $\ek\an(K_{p})_{\phi}$.
The following lemma is thus proved.

\math{Lemma.}{\isotropysubgroup}{\sl The inclusion 
$(K_{p})_{\phi}\subset K_{f}$ is a necessary and sufficient
condition for the isotropy subgroup $G_{\en}$ to be the
semidirect product 
$(K_{p})_{\phi}\timess{\rho}\ker\tau_{p}^{\ast}$.} 

\vskip0.3cm

\chapter{Submanifolds-Symplectic structure}

The coadjoint orbits of a semidirect product $G=
K\timess{\rho}V$ possess always two natural submanifolds:
the $K$-orbit $L={\script O}^{K}_{\en}$ of
$\en\an\frak g^{\ast}$ and the $V$-orbit $N={\script O}^{V}_{\en}$
of the same element. We observe that in the case where 
$K_{p}\subset K_{f}$, $\en=(f,p)\an\frak g^{\ast}$, the orbits $L$
and $Z$ are diffeomorphic: ${Z}=K/K_{p}$ 
and $L=K/(K_{p}\cap K_{f})=K/K_{p}$. 
In this section we will study the submanifolds $L$ and
$N$ as well as the symplectic structure of the coadjoint orbit 
${\script O}^{G}_{\en}$.

\math{Convention.}{\conventiontwo}{\it If $\xi$ is an element of a Lie
algebra $\frak g$, then we denote by $\xi_{\frak g^{\ast}}$ the
fundamental vector field of the coadjoint action on $\frak g^{\ast}$:
$$(\xi_{\frak g^{\ast}})_{\em}=\xi\cdot\em,\;\forall\em\an\frak
g^{\ast}.\eqn\coadjfund$$
Also, the symplectic structure we are using on a coadjoint orbit in
$\frak g^{\ast}$, is given by:
$$\omega_{\em}\big((\xi_{\frak g^{\ast}})_{\em},(\eta_{\frak
g^{\ast}})_{\em}\big)=-\em([\xi,\eta]),\eqn\coadjsymplform$$
where $\em$ is a point on the orbit.}

\vskip0.3cm

For the case of the semidirect product in which we are interested,
equation $\coadjsymplform$ takes the following form: if $\xi=(A,a)$
and $\eta=(B,b)$, $A,B\an\k$, $a,b\an V$, then:
$$\omega_{\em}\big((\xi_{\frak g^{\ast}})_{\em},(\eta_{\frak
g^{\ast}})_{\em}\big)=(A\cdot h+q\odot
a)(B)+\tau_{q}(A)(b),\eqn\vartwo$$
if $\em=(h,q)\an\frak g^{\ast}$.

We find first for $o\an L$ and $n\an N$, the tangent spaces 
$T_{o} L$ and $T_{n} N$ explicitly. Using the fact that $L$ 
and $N$ are homogeneous spaces of the groups $K$ and $V$
respectively, we have:
$$T_{o}L=\big\{\big(A\cdot h,A\cdot q\big)\;|\;A\an\k
\big\},\;o=(h,q)=(\ek,0)\cdot\en\an L\eqn\tangentspacel$$ 
and 
$$T_{n}N=\{(p\odot u,0)\;|\;u\an
V\}\cong\im\tau^{\ast}_{p}=\k_{p}^{\circ},\;n=(e,v)\cdot\en\an N.
\eqn\tangentspacen$$
We observe here that the tangent space $T_{n}N$ does not depend
on the point $n\an N$, in accordance with the affine plane
structure of $N=\{(f+p\odot v,p)\;|\;v\an V\}=
\en+\im\tau_{p}^{\ast}\times\{0\}$.

In order to obtain now a characterization for the submanifolds $L$
and $N$, we search for the orthogonal complements of the
tangent spaces $T_{o}L$ and $T_{n}N$. We find easily, using
the relations $\vartwo$ and $\tangentspacel$:
$$\big(T_{o}L\big)^{\bot}=\big\{(0,B\cdot q)\an 
T_{o}{\script O}^{G}_{\en}\;|\;B\an\k,\,B\cdot h\an q\odot V
\big\}.\eqn\orthogcompl$$
It is clear that generally, the orthogonal complement $\big(T_{o}
L)^{\bot}$ has no relation to the tangent space $T_{o}L$.
But in the case where $f$ defines a cohomology class, 
$[f]\an{\rm H}^{1}(\k,\R)$, we have $B\cdot f=0$ hence $B\cdot h=0$,
$\forall B\an\k$, which implies, by $\orthogcompl$ and $\tangentspacel$ 
that 
$$\big(T_{o}L\big)^{\bot}=T_{o}L.$$
Thus the condition $[f]\an{\rm H}^{1}(\k,\R)$ means that the submanifold
$L$ is Lagrangian.

We turn to the case of $\big(T_{n}N\big)^{\bot}$. Easy
calculation shows that 
$$\big(T_{n}N\big)^{\bot}=\big\{\big(B\cdot(f+p\odot v)+p\odot
b,0\big)\;|\;B\an\k_{p},\;b\an V\big\},$$
which clearly leads to the inclusion $T_{n}N
\subset\big(T_{n}N\big)^{\bot}$. Thus, the submanifold $N$
is always an isotropic submanifold of the coadjoint orbit ${\script
O}^{G}_{\en}$. 

When $K_{p}\subset K_{f}$, $N$ becomes a Lagrangian
submanifold. Indeed, the isomorphism $T_{n}N\cong\k^{\circ}_{p}$
implies that the dimensions of $N$ and $Z$ are equal; furthermore, if
$K_{p}\subset K_{f}$, then the tangent spaces of $N$ and $L$ at $\en$ are 
complementary, so $\dim{N}+\dim{L}=\dim{\script O}^{G}_{\en}$ and
$Z\cong L\Rightarrow\dim{Z}=\dim{L}$;
finally $2\dim{N}=\dim{\script O}_{\en}^{G}$.

We have now a useful characterization of the cotangent space 
$T^{\ast}_{q}{Z}$, $q\an Z$:

\math{Lemma.}{\naturisom}{\sl For each point $q\an Z$, the
cotangent space $T_{q}^{\ast}{Z}$ is naturally isomorphic to
the quotient $V/\ker\tau_{q}^{\ast}$.}

\undertext{\it Proof}. It follows directly from the isomorphism
$T_{q}{Z}\cong\k/\k_{q}$ and Lemma {\rangeoftau}.\QED

\vskip0.3cm

Using the previous lemma, one can investigate further the
consequences of the condition 
$K_{p}\subset K_{f}$ on the structure of the coadjoint orbit
${\script O}^{G}_{\en}$, where $\en=(f,p)\an\frak g^{\ast}$. In fact,
let $\Pi\Colon{\script O}^{G}_{\en}\rightarrow {Z}$ 
be the projection and $q\an Z$; then, by equation
$\coadjoint$, one easily finds that for $\em=(h,q)\an{\script
O}^{G}_{\en}$, the fibre $\Pi^{-1}(\Pi(\em))$ is of the form 
$$\Pi^{-1}(q)=(h,q)+q\odot V\times\{0\}.\eqn\fibre$$
Lemma {\naturisom} applied to the case $K_{p}\subset K_{f}$,
makes clear that ${\script O}_{\en}^{G}$ is isomorphic (as a manifold)
to the cotangent bundle $T^{\ast}{Z}$ when $K_{p}\subset K_{f}$.

We make now the following remark concerning the bundle $Y$ of
little-group orbits, under the condition $K_{p}\subset
K_{f}$. If $\phi=i_{p}^{\ast}f$ (notation of sections 2 
and 3), the typical fibre of $Y$ is $Y_{p}={\script
O}^{K_{p}}_{\phi}$. It is immediate that for each $\ek\an K_{p}$,
$\ek\cdot\phi=i_{p}^{\ast}(\ek\cdot f)=\phi$ which implies that 
$Y_{p}=\{\phi\}$; thus the fibre bundle $Y$ and the orbit $Z$ are
isomorphic as manifolds.

We study finally the case $[f]\an{\rm H}^{1}(\k,\R)$ and its
consequences on the structure of the coadjoint orbit ${\script
O}^{G}_{\en}$. We use the following well-known property of the
coadjoint action:

\math{Property.}{\discretorbit}{\sl If $f\an\k^{\ast}$
defines a cohomology class, $[f]\an{\rm H}^{1}(\k,\R)$, then the
coadjoint orbit ${Q}={\script O}^{K}_{f}$ is a manifold
of dimension zero.}

\vskip0.3cm

In other words, the isotropy subgroup $K_{f}$ is at the same
time open and closed in $K$; thus, using the previous discussion
and relation $\coadjoint$ we find that the fibre bundle
$\Pi\Colon{\script O}^{G}_{\en}\rightarrow{Z}$
defines a covering space of the cotangent bundle $T^{\ast}Z$. 
In particular, the bundle of little-group orbits is a covering
space of the orbit $Z$.

Let us summarize with the proposition:

\math{Proposition.}{\submanifolds}{\sl The coadjoint orbit ${\script
O}^{G}_{\en}$ of a semidirect product $G=
K\timess{\rho}V$, $\en=(f,p)\an\frak g^{\ast}$, possesses always two
natural submanifolds $L$ and $N$ which have transversal
intersection
at $\en\an L\cap N$: $L$ is the orbit of $\en$ under
the action of $K\subset G$ and $N$ the orbit of the
same point under $V\subset G$. $N$ is always an isotropic
submanifold of the coadjoint orbit ${\script O}^{G}_{\en}$.
\item{1.} Suppose that $K_{p}\subset K_{f}$; then:
\itemitem{a.} the orbit ${Z}$ is
diffeomorphic to $L$ and $N$ is a Lagrangian submanifold;
\itemitem{b.} the coadjoint orbit ${\script O}^{G}_{\en}$ is
diffeomorphic to the cotangent bundle $T^{\ast}{Z}$;
\itemitem{c.} the bundle of little-group orbits $Y$ and the orbit
$Z$ are identical.
\item{2.} Suppose that $f$ defines a cohomology class, $[f]\an{\rm
H}^{1}(\k,\R)$; then:
\itemitem{a.} $L$ and $N$ are Lagrangian submanifolds of 
${\script O}^{G}_{\en}$;
\itemitem{b.} the orbit ${\script O}^{G}_{\en}$ is a
covering space of the cotangent bundle $T^{\ast}{Z}$;
\itemitem{c.} the bundle of little-group orbits $Y$ is a covering
space of the orbit $Z$.}

\vskip0.3cm

More generally, we can define a foliation $\eusm F$ on the coadjoint
orbit ${\script O}^{G}_{\en}$ in the following way:
if $o=(h,q)=\big(\ek\cdot f,\ek\cdot p\big)\an L$,
$\ek\an K$, then we choose the leaf $\eusm F_{o}$ as $\eusm
F_{o}={\script O}_{o}^{V}=\{(h+q\odot v,q)\;|\;v\an V\}$. 
Then, using the techniques of the Proposition {\submanifolds}, one
easily proves:

\math{Proposition.}{\isotrfol}{\sl The coadjoint orbit ${\script
O}^{G}_{\en}$ of a semidirect product 
$G=K\timess{\rho}V$, $\en=(f,p)\an\frak g^{\ast}$, possesses an isotropic
foliation whose leafs are the affine spaces $\eusm 
F_{o}={\script O}^{V}_{o}$, $o\an L$ (see Proposition
{\submanifolds}). In the case where $K_{p}\subset K_{f}$,
the foliation $\eusm F$ is Lagrangian.}

\vskip0.3cm

In view of Lemma {\coadorbitstruc}, the following is immediate:

\math{Corollary.}{\corollone}{\sl The foliation $\eusm F$ of
Proposition {\isotrfol} is always regular and the quotient ${\script
O}^{G}_{\en}/\eusm F$ is equal to the bundle of little-group orbits.}

\vskip0.3cm

We observe here the analogy with the cotangent bundle: in fact, if
$T^{\ast}M$ is the cotangent bundle of a manifold $M$, then the fibres
$T^{\ast}_{m}M$, $m\an M$, define a Lagrangian foliation of
$T^{\ast}M$. We will clarify in what follows this similarity by
direct calculation of the symplectic form $\omega$ of the coadjoint
orbit ${\script O}^{G}_{\en}$ in terms of the symplectic
structures of $T^{\ast}{Z}$ and ${Q}=
{\script O}^{K}_{f}$. 

Fix now an element $\en=(f,p)\an\frak g^{\ast}$ and let 
$\sigma\Colon G\rightarrow{\script O}^{G}_{\en}$, $\sigma_{1}\Colon
G\rightarrow{Q}$ and $\sigma_{2}\Colon G\rightarrow
T^{\ast}{Z}$ be the mappings
defined by $\sigma_{1}(\ek,v)=\ek\cdot f$,
$\sigma_{2}(\ek,v)=(\ek\cdot p,[v]_{\ek\cdot p})$ and $\sigma$ is
simply the projection $G\rightarrow G/G_{\en}$;
$[v]_{\ek\cdot p}$ means the equivalence class of $v$ in the quotient
$V/\ker\tau_{\ek\cdot p}^{\ast}$, according to Lemma {\naturisom}. 

\math{Theorem.}{\semprodsymplform}{\sl The canonical symplectic
structures $\omega$, $\omega_{Q}$ and $\omega_{Z}$ of
${\script O}^{G}_{\en}$, $Q$ and $T^{\ast}Z$
respectively, are related by the following equation:
$$\sigma^{\ast}\omega=\sigma_{1}^{\ast}\omega_{Q}
+\sigma_{2}^{\ast}\omega_{Z}.\eqn\allsymplforms$$}
\indent\undertext{\it Proof}. Let
$\Omega=\sigma_{1}^{\ast}\omega_{Q}
+\sigma_{2}^{\ast}\omega_{Z}\an\Omega^{2}(G)$; $\Omega$ is a
closed 2-form. If
$\xi=(A,a),\eta=(B,b)\an\frak g$, let also $\xi^{r}$ and $\eta^{r}$ be
the corresponding right invariant vector fields on $G$. Then, easy
calculation shows that
$$\Omega(\xi^{r},\eta^{r})\;\big|_{g}=q\big(B\cdot(A\cdot v+a)\big)-
q\big(A\cdot(B\cdot v+b)\big)-h([A,B]),$$
if $q=\ek\cdot p$, $h=\ek\cdot f$ and $g=(\ek,v)$. 
On the other hand, the fact that ${\script O}^{G}_{\en}$ is the
quotient $G/G_{\en}$ implies directly
$T_{g}\sigma(\xi^{r}(g))=\xi_{\frak g^{\ast}}(\em)$, $\em=\sigma(g)$.
By Convention {\conventiontwo} and relation $\vartwo$,
this means that $\sigma^{\ast}\omega$ is exactly $\Omega$. \QED

\math{Remark.}{\noncanonical}{\it It must be emphasized that the previous
result is closely related to the choice of a point of ${\script
O}^{G}_{\en}$ (here we choose the origin $\en$
of the orbit). In this sense, the splitting of the symplectic structure
$\omega$ of the coadjoint orbit, is not canonical.}

\chapter{Coadjoint orbits, modified cotangent bundles and coisotropic
embeddings}

A cotangent bundle $T^{\ast}M$ equipped with the 2-form
$\wt{\omega}_{M}=\omega_{M}+\tau^{\ast}\alpha_{0}$, where
$\omega_{M}=d\theta_{M}$ is the canonical symplectic
form of $T^{\ast}{M}$ and $\alpha_{0}\an\Omega^{2}({M})$ is such that
$d\alpha_{0}=0$, $\tau\Colon T^{\ast}{M}\rightarrow{M}$,
is called modified
cotangent bundle. We denote the pair $(T^{\ast}{M},\omega_{M}
+\tau^{\ast}\alpha_{0})$ by $T^{\sharp}{M}$, when $\alpha_{0}$ 
is understood. 
The form $\wt{\omega}_{M}$ is always non-degenerate, as one readily
verifies, so $T^{\sharp}{M}$ is also a symplectic manifold.

This notion is closely related to physical problems. As most
important, we mention the phase space of a charged particle in general
relativity, in the presence of an external electromagnetic field, e.g.
{\sn}, {\sou} and the problem of localization of relativistic
particles of mass zero, {\det}. In this section we will see that it is
also useful in the geometry of the coadjoint orbits of a semidirect
product.

Consider an element $\en=(f,p)\an\k^{\ast}\oplus V^{\ast}$ with
$K_{p}\subset K_{f}$. Then we know that the orbit $Z$ 
forms a fibre bundle over ${Q}$ with typical fibre 
$K_{f}/K_{p}$, the orbit of $p$ under $K_{f}$. In
addition, if $\eurm p\eurm r\Colon{Z}\rightarrow{Q}$
is the projection, $\eurm p\eurm r(\ek\cdot p)=\ek\cdot f$, the 2-form
$\alpha_{0}=\eurm p\eurm r^{\ast}\omega_{Q}\an\Omega^{2}({Z})$ 
is a presymplectic structure on $Z$.

\math{Proposition.}{\modifiedcot}{\sl If the element
$\en=(f,p)\an\frak g^{\ast}$ is
such that $K_{p}\subset K_{f}$, then there exists a global
section $s\Colon{Z}\rightarrow{\script O}^{G}_{\en}$
of $\Pi\Colon{\script O}^{G}_{\en}\rightarrow{Z}$
such that $s^{\ast}\omega=\alpha_{0}$. Furthermore, if $\tau\Colon
T^{\ast}{Z}\rightarrow{Z}$ is the projection, then
there exists a symplectomorphism between 
$({\script O}^{G}_{\en},\omega)$ and $T^{\sharp}{Z}=
(T^{\ast}{Z},\omega_{Z}+\tau^{\ast}\alpha_{0})$. However,
this symplectomorphism is not canonical.}

\undertext{\it Proof}. Define $s\Colon{Z}\rightarrow{\script
O}^{G}_{\en}$ by $s(\ek\cdot p)=(\ek\cdot f,\ek\cdot p)$. The
map $s$ is well defined because if $\ek\cdot p=\el\cdot p$, then
$\el=\ek\cdot\ek\tonos$ for $\ek\tonos\an K_{p}$ and so
$s(\el\cdot p)=s(\ek\cdot p)$. Clearly, we have $\Pi\comp
s=id$, so $s$ is a section. Its tangent at $q=\ek\cdot p$ is given by
$T_{q}s(A\cdot q)=(A\cdot h,A\cdot q)=(A,0)_{\frak g^{\ast}}(h,q)$ if
$h=\ek\cdot f$; this makes clear that $s^{\ast}\omega=\alpha_{0}$ (see
{\vartwo}). 

Now by Proposition {\submanifolds} we have the diffeomorphism ${\script
O}^{G}_{\en}\cong T^{\ast}Z$. Using the notation of Theorem 
{\semprodsymplform} we observe that if $\sigma(\ek,v)=\sigma(\el,u)$, then
$\el=\ek\cdot\ek\tonos$ for $\ek\tonos\an K_{p}$ and
$v-u\an\ker\tau_{\ek\cdot p}^{\ast}$. This means that there exist two
well defined mappings $\bit I_{1}\Colon{\script O}^{G}_{\en}
\rightarrow{Q}$ and $\bit I_{2}\Colon{\script
O}^{G}_{\en}\rightarrow T^{\ast}{Z}$ satisfying $\bit
I_{i}\comp\sigma=\sigma_{i}$, $i=1,2$. Explicitly, if $\em$ is the
element of $\frak g^{\ast}$ given by relation $\coadjoint$, then
$\bit I_{1}(\em)=\ek\cdot f$ and $\bit I_{2}(\em)=(\ek\cdot
p,[v]_{\ek\cdot p})$. As a consequence, relation $\allsymplforms$
reads $\omega=\bit I_{1}^{\ast}\omega_{Q}+\bit
I_{2}^{\ast}\omega_{Z}$ because $\sigma^{\ast}$ is a
monomorphism. But $\bit I_{1}$ satisfies also $\bit I_{1}=\eurm p\eurm
r\comp\tau\comp\bit I_{2}$ which gives $\omega=\bit
I_{2}^{\ast}(\omega_{Z}+\tau^{\ast}\alpha_{0})$. Finally, it
is elementary to verify that $\bit I_{2}$ is a diffeomorphism between
${\script O}^{G}_{\en}$ and $T^{\ast}{Z}$; this means that we
have a symplectomorphism between ${\script O}^{G}_{\en}$ and
$T^{\sharp}{Z}$ which is not canonical
because it depends on the choice of a point $\en=(f,p)$ of the
coadjoint orbit in question, for which the condition $K_{p}\subset
K_{f}$ is satisfied. \QED

\vskip0.3cm

We summarize our knowledge on the structure of ${\script O}^{G}_{\en}$
when $K_{p}\subset K_{f}$ in the commutative
diagram:

\vskip0.5cm

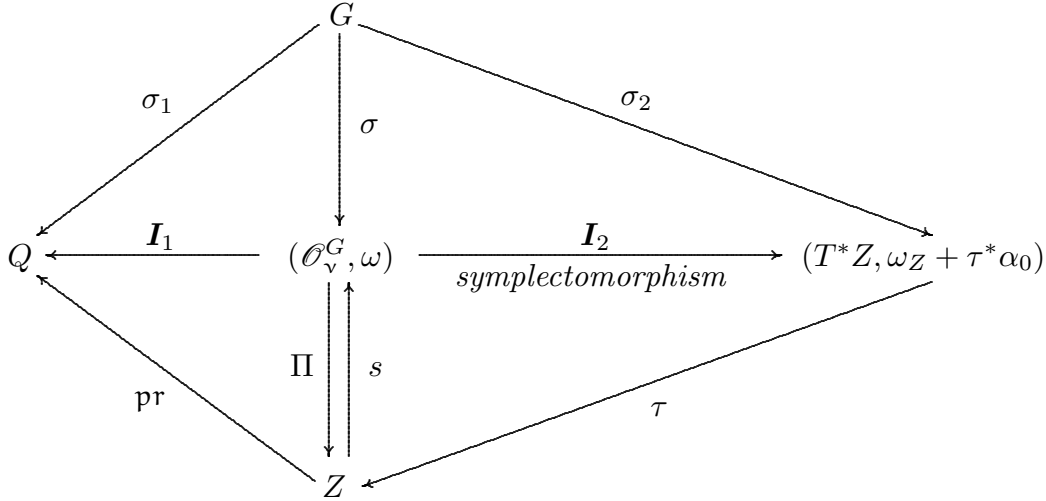
\centerline{{\beginpicture
\setcoordinatesystem units <0.70000cm,0.70000cm>
\setplotsymbol ({\fiverm .})
\setlinear

\arrow <1mm> [0.1,1] from 11.589 16.669 to 17.304 16.669
\arrow <1mm> [0.1,1] from 9.049 16.669 to  5.7 16.669
\arrow <1mm> [0.1,1] from 10.001 20.320 to 5.556 16.986
\arrow <1mm> [0.1,1] from 10.636 20.320 to 19.685 16.986
\arrow <1mm> [0.1,1] from 10.33 20.161 to 10.33 17.15 
\arrow <1mm> [0.1,1] from 10.160 16.24 to 10.160 13.494
\arrow <1mm> [0.1,1] from 10.478 13.494 to 10.478 16.24
\arrow <1mm> [0.1,1] from 19.685 16.251 to 10.7 13.018 
\arrow <1mm> [0.1,1] from 9.94 13.1 to 5.556 16.351

\put{$G$} [lB] at 10.160 20.320
\put{$(T^{\ast}{Z},\omega_{Z}+\tau^{\ast}\alpha_{0})$} [lB]
at 17.621 16.510
\put{$Q$} [lB] at  5.080 16.510
\put{$({\script O}^{G}_{\en},\omega)$} [lB] at  9.5 16.510
\put{$Z$} [lB] at 10.05 12.859
\put{$\tau$} [lB] at 15.240 14.129
\put{$\sigma$} [lB] at 10.65 18.6
\put{\it symplectomorphism} [lB] at 12.15 16.192
\put{$\Pi$} [lB] at  9.55 14.764
\put{$s$} [lB] at 10.78 14.764
\put{$\sigma_{1}$} [lB] at  7.2 19
\put{$\sigma_{2}$} [lB] at 14.764 19.050
\put{$\bit I_{2}$} [lB] at 14.129 16.828
\put{$\bit I_{1}$} [lB] at  7.26 16.828
\put{$\eurm p\eurm r$} [lB] at  7.1 14.3
\endpicture}
}

\vskip0.2cm

But there exist additional properties of the coadjoint orbit ${\script
O}^{G}_{\en}$ when $K_{p}\subset K_{f}$. More
precisely:

\math{Theorem.}{\coisotropicemb}{\sl Let 
$G=K\timess{\rho}V$ be a semidirect product and $\en=(f,p)\an\frak
g^{\ast}$ such that $K_{p}\subset K_{f}$. Then the
reduction of ${\script O}^{G}_{\en}$ by the submanifold $L$ of
Proposition {\submanifolds} is symplectomorphic to the coadjoint orbit
${Q}=K/K_{f}$. Furthermore, the dual $E^{\ast}$ of
the characteristic distribution on $L$ defines a symplectic manifold
which is a vector bundle over $Z$ and the zero section
$s_{0}\Colon{Z}\rightarrow E^{\ast}$ is a coisotropic
embedding of the presymplectic manifold $Z$.}

\undertext{\it Proof}. The section $s$ is a diffeomorphism between
$Z$ and $L$ and therefore defines a pre-symplectomorphism
$s\Colon({Z},\alpha_{0})\rightarrow(L,i^{\ast}_{L}\omega)$,
$i_{L}\Colon L\rightarrow{\script O}^{G}_{\en}$ is the
inclusion, because $s^{\ast}i_{L}^{\ast}\omega=(i_{L}\comp
s)^{\ast}\omega=\alpha_{0}$. This ensures that the characteristic
distribution $E$ of $i^{\ast}_{L}\omega$ on $L$, $E=TL\cap(TL)^{\bot}$
is isomorphic to the
kernel of $\alpha_{0}$ over $Z$. But if we reduce $Z$ by
$\ker\alpha_{0}$ we obtain the coadjoint orbit $K/K_{f}$. Thus,
reduction of ${\script O}^{G}_{\en}$ by $L$ gives the symplectic
manifold $Q$.

For the case we are studying, one easily finds that if $o=(\ek\cdot
f,\ek\cdot p)=(h,q)\an L$, then, see $\orthogcompl$,
$T_{o}L\cap(T_{o}L)^{\bot}=\{(0,A\cdot
q)\;|\;A\an\k_{h}\}$. Thus, we conclude that the dual $E^{\ast}$ of
the characteristic distribution $E$, is a vector bundle over $Z$ 
with typical fibre $(\k_{f}/\k_{p})^{\ast}$. Consider now a point
$x\an E$; let $q=\ek\cdot p$ be its projection on $Z$ and
$h=\ek\cdot f=\eurm p\eurm r(q)$. Then, by construction of $E^{\ast}$,
the tangent space $T_{x}E^{\ast}$ admits the decomposition:
$T_{x}E^{\ast}=T_{h}{Q}
\oplus(\k_{h}/\k_{q})\oplus(\k_{h}/\k_{q})^{\ast}$. This makes clear
that there exists a symplectic structure $\omega_{E^{\ast}}$ on
$E^{\ast}$ given by
$(\omega_{E^{\ast}})_{x}(\xi_{1}+\alpha_{1}+\beta_{1},
\xi_{2}+\alpha_{2}+\beta_{2})=(\omega_{Q})_{h}
(\xi_{1},\xi_{2})+\beta_{1}(\alpha_{2})-\beta_{2}(\alpha_{1})$,
$\xi_{i}\an T_{h}{Q}$, $\alpha_{i}\an(\k_{h}/\k_{q})$,
$\beta_{i}\an(\k_{h}/\k_{q})^{\ast}$, $i=1,2$.
Let now $s_{0}\Colon{Z}\rightarrow E^{\ast}$ be the zero
section and ${\ygot Z}=s_{0}({Z})$. Then each tangent vector
$v\an T_{x}{\ygot Z}$ admits the decomposition $v=\eta+\gamma+0$,
$\eta\an T_{h}{Q}$, $\gamma\an(\k_{h}/\k_{q})$. Consequently,
the orthogonal complement of this tangent space will be given by 
$(T_{x}{\ygot Z})^{\bot}=\{0+\alpha+0\an
T_{x}E^{\ast}\;|\;\alpha\an(\k_{h}/\k_{q})\}\subset T_{x}{\ygot Z}$, 
which completes the proof of the theorem. \QED

\vskip0.3cm

We make some comments on the
somewhat peculiar condition $K_{p}\subset K_{f}$. In this way stated,
this condition depends on a point $(f,p)$ of the coadjoint orbit, but
as one easily verifies, it implies that
$$\ek\cdot h-h\an\im\tau_{q}^{\ast},\;\forall\ek\an K_{q},
\eqn\varisotropy$$
for an arbitrary point $(h,q)\an{\script O}^{G}_{\en}$. Now, condition
$\varisotropy$ is equivalent to saying that all the little-group orbits
are trivial because $i_{q}^{\ast}(\ek\cdot h)=i_{q}^{\ast}(h)$,
$\forall\ek\an K_{q}$. Otherwise stated, condition $K_{p}\subset
K_{f}$ implies that the orbit $Z$ coincides with the bundle of
little-group orbits $Y$. 

The converse now is not in general true, that is the condition
$${Z}=Y$$
does not in general implies the existence of an
element $(f,p)\an{\script O}^{G}_{\en}$ such that $K_{p}\subset
K_{f}$. Let us explain why. If ${Z}=Y$ then
condition $\varisotropy$ is valid, so given $(h,q)\an{\script
O}^{G}_{\en}$ and $\ek\an K_{q}$, we have $\ek\cdot h-h=q\odot
v(\ek)$, where $v\Colon K_{q}\rightarrow V$ is such that its
equivalence class $[v]\Colon K_{q}\rightarrow V/\ker\tau_{q}^{\ast}$
belongs to $Z^{1}(K_{q},V/\ker\tau_{q}^{\ast})$. This is evident by
direct calculation using the fact that $\ker\tau_{q}^{\ast}$ is
$K_{q}$-invariant, which induces a representation $K_{q}\rightarrow
GL(V/\ker\tau_{q}^{\ast})$. If now there exists
an element $(h,q)\an{\script O}^{G}_{\en}$ such that
${\rm H}^{1}(K_{q},V/\ker\tau_{q}^{\ast})=0$,
then we can find an element $(f,p)$
of the same orbit with the property $K_{p}\subset K_{f}$. In fact, in
this case we have always $q\odot v(\ek)=q\odot(\ek\cdot v_{0}-v_{0})$
for a fixed element $v_{0}\an V$. Choosing thus $p=q$ and $f=h-q\odot
v_{0}$, we have the desired result.

\chapter{Pukanszky's condition and the semidirect product}

Let us first state some definitions and results {\det} about polarizations
and Pukanszky's condition that will be used in the sequel.

Let $G$ be a Lie group, $\frak g$ its Lie algebra and $\en$ 
an  element of $\frak g^{\ast}$. 
Given a subspace $\frak a\subset\frak g$ which contains
the Lie algebra $\frak g_{\en}$ of the isotropy subgroup
$G_{\en}$ with respect to the coadjoint action,
we define the symplectic orthogonal $\frak
a^{\perp}$ by $$\frak a^{\perp}=\{X\an\frak g\;|\;\en([X,Y])=0,\;
\forall Y\an\frak a\}.\eqn\symplorthogonal$$
If we note by $\frak g^{\complex}$
the complexification of $\frak g$ and by $\frak g^{\complex}\na
\em\mapsto\bar{\em}\an\frak g^{\complex}$ the complex 
conjugation, we may extend this notion immediately for subspaces
of $\frak g^{\complex}$ which contain $\frak g^{\complex}_{\en}$.  

We say now that the complex Lie subalgebra $\frak h$ of $\frak
g^{\complex}$ is a polarization with respect to $\en\an\frak g^{\ast}$
if $\frak h$ contains $\frak g^{\complex}_{\en}$, is invariant under
the adjoint action of $G_{\en}$, $\frak h^{\perp}=\frak h$ and
$\frak h+\bar{\frak h}$ is a Lie subalgebra of $\frak g^{\complex}$.
Each algebraic polarization $\frak h$ corresponds to a
$G$-invariant geometric polarization $\eusm F$, the correspondence
being given by $\frak h\cdot\en=\eusm
F_{\en}\subset(T_{\en}{\script O}_{\en}^{G})^
{\complex}$. The condition on the symplectic orthogonal $\frak h^{\perp}$
can be restated as follows:
$$(\frak h^{\bot}=\frak h)\Longleftrightarrow
\big(\dim_{\complex}\frak h={1\over{2}}(\dim\frak g+\dim\frak g_{\en})
\quad\hbox{and}\quad
\en([\,\frak h,\frak h\,])=0\big).\eqn\isotrdimmax$$
\indent To each polarization $\frak h$ we can associate two real
Lie subalgebras $\don\subset\frak e$ of $\frak g$ defined by 
$$\don=\frak h\cap\frak g\quad\hbox{and}\quad\frak 
e=(\,\frak h+\bar{\frak h}\,)\cap\frak g.\eqn\realsubalg$$
We denote also by $D_{0}\subset E_{0}$ the connected Lie
subgroups of $G$ whose Lie algebras are $\don$ and $\frak e$
respectively. The conditions on the polarization $\frak h$ ensure that
the subsets $D=D_{0}\cdot G_{\en}\subset E=E_{0}\cdot
G_{\en}$ are subgroups of $G$.

\math{Lemma, {\kost}.}{\polarisationlemma}{\sl (1) 
$\don^{\perp}=\frak e$; (2) the groups
$D$ and $D_{0}$ are closed Lie subgroups
of $G$ with the same Lie algebra $\don$; (3) $i^{\ast}_{\don}\en$ is
invariant under the coadjoint action of $D$;
(4) $\frak e^{\circ}$ and $\en+\frak e^{\circ}$
are invariant under the coadjoint action of the subgroup $D$;
(5) if $E$ is a Lie subgroup of $G$, 
then its Lie algebra is $\frak e$.}

\vskip0.3truecm

Pukanszky's condition now, is a supplementary condition on the
polarization $\frak h$. The following lemma gives three equivalent 
variants of this condition. 

\math{Lemma {\bfsl (Pukanszky's condition)}, {\det}.}{\pukcond}
{\sl The following conditions are equivalent:
(1) $\en+\frak e^{\circ}
\subset{\script O}_{\en}^{G}$;
(2) $D\cdot\en=\en+\frak e^{\circ}$
(3) $D\cdot\en$ is closed in $\frak g^{\ast}$.}

\vskip0.3cm

Consider now the case where the Lie group $G$ is a semidirect
product, $G=K\timess{\rho}V$ (notation of section 2).
Then, its Lie algebra is $\frak g=\k\opluss{\rho}V$ and the
corresponding complexified Lie algebra $\frak g^{\complex}={\ygot
k}^{\complex}\opluss{\rho} V^{\complex}$. We are interested in
polarizations of $\frak g^{\complex}$ (with respect to $\en\an\frak
g^{\ast}$) which are of the form $\frak h=\frak
a\opluss{\rho}V^{\complex}$, $\frak a\subset\k^{\complex}$. Although
this type of polarization seems to be very special, it leads to quite
interesting results as we shall see in the sequel.

We examine first the restrictions imposed to the subspace $\frak a$
by the fact that $\frak h$ is a polarization. We find
successively: 
\item{(1)} {\it $\frak h$ is a subalgebra of $\frak g^{\complex}$}.
Then, $[\,\frak h,\frak h\,]\subset\frak h$ which implies 
$[\,{\frak a},{\frak a}\,]\opluss{\rho}[\,{\frak a},V^{\complex}\,]\subset
{\frak a}\opluss{\rho}V^{\complex}$. Thus, $\frak a$ must be a Lie 
subalgebra of $\k^{\complex}$.
\item{(2)} {\it $\frak h^{\bot}=\frak h$}. Equivalently, we have the
relation $\isotrdimmax$. By direct calculation of the dimensions
appearing in $\isotrdimmax$, we obtain:
$\dim_{\complex}\frak h=\dim_{\complex}{\frak a}+\dim V$ and 
$\dim\frak g+\dim\frak g_{\en}=2\dim
V+\dim\k_{p}+\dim(\k_{p})_{\phi}$. Thus,
$$\dim_{\complex}{\frak a}
={1\over{2}}\big(\dim\k_{p}+\dim(\k_{p})_{\phi}\big).\eqn\dimeneq$$
We have one more restriction coming from the condition
$\en([\,\frak h,\,\frak h])=0$. Indeed, using relations $\coadjsymplform$
and $\vartwo$,
this condition gives: $A\cdot p=0$
and $A\cdot f+p\odot a\an\frak a^{\circ}$, $\forall A\an\frak a, a\an
V^{\complex}$. But if this is the case, we have $\frak
a\subset\k^{\complex}_{p}$ and since $p\odot
a\an(\im\tau_{p}^{\ast})^{\complex}=(\k_{p}^{\complex})^{\circ}$ 
(Lemma {\rangeoftau}), for each $a\an V^{\complex}$, we obtain
necessarily $A\cdot f\an\frak a^{\circ}$, $\forall A\an\frak a$, or,
$f([\,\frak a,\frak a\,])=0$. Taking into account relation {\dimeneq}
as well as the facts $\frak a\subset\k_{p}^{\complex}$ and $f([\,\frak
a,\frak a\,])=0$, we conclude that $\frak a^{\bot}=\frak a$ (the
symplectic orthogonal being taken with respect to
$\phi=i_{p}^{\ast}f$). 
\item{(3)} {\it $\frak h$ is invariant under the adjoint action of
$G_{\en}$}. By equation $\adjoint$ it is immediate that $\frak
a$ must be invariant under the adjoint action of $(K_{p})_{\phi}$.
\item{(4)} {\it $\frak h+\bar{\frak h}$ is a
Lie subalgebra of $\frak g^{\complex}$}.
Equivalently, it suffices to
demand $[\,\frak h,\bar{\frak h}\,]\subset\frak
h+\bar{\frak h}$, because $\frak h$ is a Lie subalgebra of $\frak
g^{\complex}$. But the last bracket is equal to $[\,\frak
a,\bar{\frak a}\,]\oplus[\,\frak a+\bar{\frak a},V^{\complex}]$; 
thus we require only $[\,{\frak a},\bar{\frak
a}\,]\subset{\frak a}+\bar{\frak a}$, because the last two terms 
belong already in $V^{\complex}$.

\vskip0.3cm

Conversely, suppose that $\frak a$ is a polarization of
$\k_{p}^{\complex}$ (with respect to $\phi=i_{p}^{\ast}f$). In that
case, it is easy to reverse the previous reasonings and deduce that 
$\frak h={\frak a}\opluss{\rho} V^{\complex}$ is a polarization (with
respect to $\en$) for the Lie algebra $\frak g^{\complex}$ (see also
{\ra}).

\math{Proposition.}{\semprodpol}{\sl Let $G=K\timess{\rho}
V$ be a semidirect product, $\en=(f,p)\an\break \frak g^{\ast}
=\k^{\ast}\opluss{\rho} V^{\ast}$ an element of the dual of its Lie
algebra and $\frak g_{\en}$ the isotropy subalgebra of $\en$ with
respect to the coadjoint action. Then, a subspace $\frak
h={\frak a}\opluss{\rho}V^{\complex}\subset\frak g^{\complex}$ is a
polarization for the Lie algebra $\frak g^{\complex}$ with respect to
$\en$ if and only if $\frak a$ is a polarization of
$\k_{p}^{\complex}$ with respect to $\phi=i_{p}^{\ast}f$.}

\vskip0.3cm

Next, we show that an analogous phenomenon occurs with
respect to the validity of Pukanszky's condition. More precisely, the
validity of Pukanszky's condition for polarizations of a semidirect
product $G=K\timess{\rho} V$ described in Proposition {\semprodpol},
reduces to validity of the same condition for polarizations of 
the isotropy Lie subalgebra $\k_{p}$.
 
Suppose that the polarization $\frak h={\frak a}\opluss{\rho}
V^{\complex}$ of $\frak g^{\complex}$ satisfies Pukanszky's condition.
For the case of the semidirect product we are
interested in, the (real) subalgebras $\don$ and $\frak e$ (see equation
$\realsubalg$) will be:
$\don=\frak h\cap\frak g={\frak a}\cap\k\opluss{\rho} V=\frak
p\opluss{\rho} V$ and
$\frak e=(\,\frak h+\bar{\frak h}\,)\cap\frak g=(\,{\frak
a}+\bar{\frak a}\,)\cap\k\opluss{\rho}
V=\frak q\opluss{\rho} V$. The connected Lie subgroup $D_{0}$
whose Lie algebra is equal to $\don$, is a semidirect product 
$D_{0}=P_{0}\timess{\rho}V$, where $P_{0}$ is the closed,
connected and simply connected
subgroup of $K_{p}$ whose Lie algebra is $\frak p$.
The validity of Pukanszky's condition for the polarization $\frak h$
implies the relation $D_{0}\cdot\en=\en+\frak e^{\circ}$ (see
Lemma {\pukcond}(2)). Consider then the element $d=(\el,v)\an
D_{0}$ and let $\en=(f,p)$ as previously: $d\cdot\en
=\big(\el\cdot f+\el\cdot p\odot v,\el\cdot p\big)$.
We know that the difference $d\cdot\en-\en$ must be contained in 
$\frak e^{\circ}$, so 
$$\big(\el\cdot f-f+\el\cdot p\odot v\big)\an\frak
q^{\circ}\eqn\nulcondone$$ and
$$\big(\el\cdot p-p\big)\an V^{\circ}=0.\eqn\nulcondtwo$$
It follows that $\el\cdot p=p\Rightarrow\el\an K_{p}$
which is already satisfied because $\k_{p}\supset\frak p={\frak a}
\cap\k\supset(\k_{p})_{\phi}$ and $P_{0}$ is connected.
Furthermore $\big(\el\cdot f-f+p\odot v\big)\an\frak q^{\circ}$; but 
$\frak p\subset\frak q\subset\k_{p}\Rightarrow\frak p^{\circ}\supset\frak
q^{\circ}\supset\k_{p}^{\circ}$ and we obtain $p\odot
v\an\im\tau_{p}^{\ast}=\k_{p}^{\circ}\subset\frak q^{\circ}$. Thus, we
must have $(\el\cdot f-f)\an\frak q^{\circ}$, or
$(\el\cdot\phi-\phi)\an i_{p}^{\ast}\frak q^{\circ}$.

As a result, the condition $D_{0}\cdot\en=\en+\frak e^{\circ}$
(Pukanszky's condition) implies $P_{0}\cdot\phi
=\phi+i_{p}^{\ast}\frak q^{\circ}$ which is exactly Pukanszky's condition
for the polarization $\frak a$ of $\k_{p}^{\complex}$ because
$i_{p}^{\ast}\frak q^{\circ}$ is the annihilator of $\frak q$ in
$\k_{p}$. Finally, the previous analysis shows easily that the
converse is also true, that is if the polarization $\frak a$ of
$\k_{p}^{\complex}$ satisfies Pukanszky's condition, then the same is
true for the polarization $\frak h={\frak a}\opluss{\rho}
V^{\complex}$ of $\frak g^{\complex}$. We have thus proved:

\math{Theorem.}{\semprodpuk}{\sl Let $\frak h={\frak a}\opluss{\rho}
V^{\complex}$ be a polarization of the complexified Lie algebra $\frak
g^{\complex}$ of a semidirect product $G=K\timess{\rho}
V$ with respect to an element $\en=(f,p)\an\frak g^{\ast}$;
equivalently 
$\frak a$ is a polarization of $\k_{p}^{\complex}$ with respect to
$\phi=i_{p}^{\ast}f$. Then, $\frak h$ satisfies Pukanszky's condition
if and only if $\frak a$ satisfies it as well.}

\math{Scholium.}{\schol}{\it The Lie algebra $\frak g$ of a
semidirect product is a special case of extension of a Lie algebra
$\k$ by an abelian Lie algebra $V$. 
The complex Lie subalgebra $\frak h$ studied in this section
is also of this special type. But one can
reconsider Proposition {\semprodpol} and Theorem {\semprodpuk} in the
following way. Let $\frak
g^{\complex}=\k^{\complex}\opluss{\rho}V^{\complex}$ and $\frak
h\subset\frak g^{\complex}$ be a Lie subalgebra. If $\frak a$ is the
image of the homomorphism $\frak g^{\complex}\rightarrow\k^{\complex}$ 
restricted to $\frak h$, then $\frak
a$ is a Lie subalgebra of $\k^{\complex}$, 
see $\commut$. In other words,
$\frak h$ is an extension of $\frak a$ by a vector subspace of
$V^{\complex}$, the kernel of $\frak h\rightarrow\frak a$. Considering
now Lie subalgebras $\frak h$ which are extensions of Lie subalgebras
of $\k^{\complex}$ by $V^{\complex}$,
$$0\longrightarrow V^{\complex}\longrightarrow\frak h
\longrightarrow\frak a\longrightarrow 0,\eqn\Hexactseq$$
one easily verifies that Proposition {\semprodpol} and Theorem
{\semprodpuk} remain still valid if we replace simply $\frak h=\frak
a\opluss{\rho}V^{\complex}$ by the exact sequence $\Hexactseq$. Finally,
we note that the corresponding geometric polarization belongs to
the category of polarizations studied in {\rm{\ra}}.} 

\vskip0.3cm

As an illustration, we then construct explicitly, for an
arbitrary semidirect product, a polarization ``trivial" in some sense,
which satisfies Pukanszky's condition 
(making of course appropriate choices). Thanks to Proposition
{\semprodpol} and Theorem {\semprodpuk}, it is sufficient to construct
a polarization $\frak a$ of $\k_{p}$ satisfying the same condition.

Suppose that $[f]\an{\rm H}^{1}(\k_{p},\R)$ and let ${\frak 
a}=\bar{\frak a}=\k_{p}^{\complex}$. Clearly, $\frak a$ is a real
subalgebra of $\k_{p}^{\complex}$, invariant under the adjoint action
of the isotropy subgroup $(K_{p})_{\phi}$ which in this case, is
the union of connected components of $K_{p}$. 
Furthermore, $2\dim_{\complex}{\frak
a}=2\dim\k_{p}=\dim\k_{p}+\dim(\k_{p})_{\phi}$, for
$(\k_{p})_{\phi}=\{A\an\k_{p}\;|\;A\cdot\phi=0\}=\k_{p}$ and the Lie 
subalgebras $\frak p$ and $\frak q$ coincide: $\frak p=\frak
q=\k_{p}$. Finally, the element $\phi=i_{p}^{\ast}f$ vanishes on all
the brackets $[A,B]$, for each $A,B\an\k_{p}$, because we have always
$[f]\an{\rm H}^{1}(\k_{p},\R)$. We deduce that ${\frak
a}=\k_{p}^{\complex}$ is a (real) polarization of $\k_{p}$.

Next, for each $\el\an P_{0}=(K_{p})_{0}$ we find:
$\el\cdot\phi=\phi$, so $\el\cdot\phi=\phi+i_{p}^{\ast}\frak 
q^{\circ}$ (here $i_{p}^{\ast}\frak
q^{\circ}=i_{p}^{\ast}\k_{p}^{\circ}=0$
is the annihilator of $\k_{p}$ in $\k_{p}^{\ast}$). 

\math{Corollary.}{\realpol}{\sl For each semidirect product
$G=K\timess{\rho}V$, the
coadjoint orbit of the element $\en=(f,p)\an\frak g^{\ast}$ 
with $[f]\an{\rm H}^{1}(\k_{p},\R)$, admits a real polarization
satisfying Pukanszky's condition.}

\chapter{Symplectic induction}

We recall here the method of symplectic induction {\kazkost} 
which will be very important to a deeper geometrical investigation
of the coadjoint orbits of semidirect products. More on
this method can be found in {\det}. Notice that the induction of
Hamiltonian actions appears independently in {\zak}.

Let $(M,\omega)$ be a symplectic manifold and $H$ a closed Lie
subgroup of a Lie group $G$. Suppose we have a (left)
Hamiltonian action $\Phi\Colon H\times M\rightarrow M$ which
admits an equivariant momentum map $J_{M}\Colon
M\rightarrow\frak h^{\ast}$, where $\frak h$ is the Lie algebra of 
$H$. The aim of the symplectic induction method is to construct a
symplectic manifold, denoted as $M_{ind}$, 
on which the group $G$ acts in a Hamiltonian way with
equivariant momentum map $J_{ind}\Colon M_{ind}\rightarrow \frak
g^{\ast}$, where $\frak g$ is the Lie algebra of $G$. 

In order to construct the Hamiltonian space $(M_{ind},\omega_{ind},
G,J_{ind})$, we proceed as follows. Using the natural isomorphism
$T^{\ast}G\cong G\times\frak g^{\ast}$ (obtained by
identifying $\frak g^{\ast}$ with the left-invariant 1-forms on 
$G$), we obtain a left action
$\check{\Phi}$ of $H$ on $\check{M}=M\times T^{\ast}G$
given by:
$$\check{\Phi}_{h}(m,g,\em)=(\Phi_{h}(m),gh^{-1},
h\cdot\em),\;\forall h\an H,(m,g,\em)\an M\times
T^{\ast}G.\eqn\combaction$$
This action is symplectic for the symplectic structure 
$\check{\omega}=\pi_{1}^{\ast}\omega_{M}+\pi_{2}^{\ast}d\theta_{G}$ 
on $\check{M}$ if $\pi_{1}\Colon\check{M}\rightarrow M$ and
$\pi_{2}\Colon\check{M}\rightarrow T^{\ast}G$ 
are the projections:
$\check{\Phi}_{h}^{\ast}\check{\omega}=\check\omega$;
$\check{\Phi}$ is also
proper because $H$ is closed. Furthermore, it 
admits an equivariant momentum map $J_{\check{M}}\Colon
\check{M}\rightarrow\frak h^{\ast}$ equal to $J_{\check{M}}
=\pi_{1}^{\ast}J_{M}+\pi_{2}^{\ast}J_{H}$, where
$J_{H}$ is the momentum map for the cotangent lift of
the right action of
$H$ on $G$. If $i_{\frak h}\Colon\frak h\hookrightarrow\frak g$
is the inclusion and $i_{\frak h}^{\ast}\Colon\frak
g^{\ast}\rightarrow\frak h^{\ast}$ the corresponding projection, then
this momentum map is given by
$$J_{H}(g,\em)=-i_{\frak h}^{\ast}\em.$$
The element $0\an\frak h^{\ast}$ is a regular
value for the momentum map $J_{\check{M}}$ and so the quotient 
$M_{ind}=J_{\check{M}}^{-1}(0)/H$ will be a symplectic manifold 
(Marsden-Weinstein reduction). We call $M_{ind}$ induced symplectic
manifold and we denote it by ${\rm Ind}^{G}_{H}M$;
$\omega_{ind}$ will denote the symplectic structure of $M_{ind}$.

In order to obtain a Hamiltonian action of $G$ on $M_{ind}$ 
we let the group $G$ act trivially on $M$; we consider also the
canonical lift to $T^{\ast}G$ of the left multiplication on $G$.
Then we have a Hamiltonian action of $G$ on $\check{M}$ 
with equivariant momentum map $\check{J}\Colon{\check{M}}
\rightarrow \frak g^{\ast}$ given by 
$$\check{J}(m,g,\em)=g\cdot\em.\eqn\indmomentmap$$
This action commutes with the action of $H$ on $\check{M}$ and
leaves invariant the momentum map $J_{\check{M}}$, so a symplectic
action of $G$ is induced on $M_{ind}$. Since 
the momentum map $\check{J}$ is $H$-invariant, 
it descends as an equivariant momentum map $J_{ind}
\Colon M_{ind}\rightarrow\frak g^{\ast}$  for the action of $G$ 
on $M_{ind}={\rm Ind}^{G}_{H}M$. 

\math{Proposition, {\det}.}{\indmanifold}{\sl 
The induced symplectic manifold 
$M_{ind}={\rm Ind}^{G}_{H}M$ is a fibre bundle over 
$T^{\ast}(G/H)$ with typical fibre the symplectic
manifold $M$. Moreover, the restriction of $\omega_{ind}$ 
to a fibre yields the original symplectic structure
$\omega_{M}$ on $M$.}

\vskip0.3cm

Let us note here that the symplectic induction procedure can be
carried out without using the trivialization of $T^{\ast}G$, see
{\land}.

If we perform now the symplectic induction for $M={\rm point}$, then
the induced symplectic manifold is isomorphic as a manifold to 
$T^{\ast}(G/H)$; we can extend the isomorphism to the
symplectic category if we modify the natural symplectic structure
$d\theta_{G/H}$ of $T^{\ast}(G/H)$ by a 
``magnetic" term, that is by the pull-back of an appropriate closed 
2-form $\beta_{0}\an\Omega^{2}(G/H)$, {\det}. Thus, the symplectic
induction over a point leads to the modified cotangent bundle
$(T^{\ast}(G/H),d\theta_{G/H}+\tau^{\ast}\beta_{0})$, where 
$\tau\Colon T^{\ast}(G/H)\rightarrow G/H$ is the cotangent projection.

\chapter{The structure of coadjoint orbits}

We will use now the results of {\det} on the structure of the coadjoint
orbits endowed with a polarization satisfying Pukanszky's condition in
order to analyze further the structure of the semidirect product
coadjoint orbits.

Let us recall first the principal results of this article. 
If $G$ is a Lie
group, $\en\an\frak g^{\ast}$ an element of the dual of its Lie
algebra and $\frak h$ a polarization with respect to $\en$, we have
two real Lie subalgebras $\don$ and $\frak e$ of $\frak g$ canonically
associated to $\frak h$ (see relations $\realsubalg$). If we define
$\eurm F=\big(G\times(\en+\frak e^{\circ})\big)/D$ (for
the property of the group $D$ see Lemma {\polarisationlemma}), then
$\eurm F$ is a vector subbundle of $T^{\sharp}(G/D)$ with symplectic
form $\omega_{\eurm F}$ obtained by restriction of the symplectic form
of $T^{\sharp}(G/D)$; further, there is a
symplectic action of $G$ on $\eurm F$ admitting a momentum map
$J_{\eurm F}\Colon\eurm F\rightarrow\frak g^{\ast}$ which may be
calculated via $\indmomentmap$.

\math{Proposition.}{\coadjpukan}{\sl The following four conditions on
the polarization are equivalent: 
\item{(1)} Pukanszky's condition.
\item{(2)} The momentum map $J_{\eurm F}\Colon\eurm
F\rightarrow\frak g^{\ast}$ is onto ${\script O}_{\en}^{G}$.
\item{(3)} The symplectic action of $G$ on $\eurm F$ is
transitive.
\item{(4)} $J_{\eurm F}$ is a symplectomorphism between $(\eurm
F,\omega_{\eurm F})$ and ${\script O}_{\en}^{G}$.}

\math{Proposition.}{\coadjpukansuppl}{\sl If $E$ is a closed
subgroup of $G$ and if Pukanszky's condition is satisfied, then
there exists a commutative diagram 

\vskip0.5cm	
\centerline{\beginpicture
\setcoordinatesystem units < 0.800cm, 0.800cm>
\setplotsymbol ({\fiverm .})
\setlinear
\setplotsymbol ({\fiverm .})

\arrow <1mm> [0.1,1] from 10.763 20.352 to 14.573 20.352
\arrow <1mm> [0.1,1] from 9.493 20.352 to  5.683 20.352
\arrow <1mm> [0.1,1] from 10.763 16.701 to 14.8 16.701
\arrow <1mm> [0.1,1] from 15.526 19.876 to 15.526 17.177
\arrow <1mm> [0.1,1] from 10.128 19.876 to 10.128 17.177
\arrow <1mm> [0.1,1] from 4.731 19.876 to  4.731 17.177
\putrule from  5.366 16.796 to  9.493 16.796
\putrule from  5.366 16.669 to  9.493 16.669

\put {${\script O}_{\en}^{G}$} [lB] at  9.81 20.193
\put {$T^{\ast}(G/D)$} [lB] at  3.8 20.193
\put {$T^{\ast}(G/E)$} [lB] at 14.732 20.193
\put {$G/E$} [lB] at 15.15 16.6
\put {$G/D$} [lB] at  4.27 16.6
\put {$G/D$} [lB] at  9.72 16.6
\put {$i_{\en}$} [lB] at  7.557 20.5
\put {$\epi_{E}$} [lB] at 12.351 20.5
\put {$\epi_{\en}$} [lB] at 10.382 18.542
\put {$\tau_{D}$} [lB] at  4.985 18.510
\put {$\tau_{E}$} [lB] at 15.748 18.542
\endpicture

}
\vskip-1.6cm

with the following properties:
\item{(1)} $(i_{\en},\epi_{\en})$ is the identification of the coadjoint
orbit ${\script O}^{G}_{\en}$ as a symplectic subbundle of
$T^{\ast}(G/D)$ according to Proposition {\coadjpukan}. 
\item{(2)} $\epi_{E}\Colon{\script O}^{G}_{\en}
\rightarrow T^{\ast}(G/E)$ is a fibre bundle whose fibres
together with the restricted symplectic form, are symplectomorphic to
the (pseudo-)K\"ahler space $E/D$.}

\vskip0.3cm

Consider now the case where the Lie group $G$ is a semidirect
product, $G=K\timess{\rho}V$ and let $\frak h$ be a
polarization of $\frak g^{\ast}$ with respect to $\en=(f,p)\an\frak
g^{\ast}$. We are always interested in polarizations of the form
$\frak h=\frak a\opluss{\rho}V^{\complex}$ satisfying Pukanszky's
condition. By Proposition {\semprodpol}, $\frak a$ is a polarization
of $\k_{p}^{\complex}$ with respect to $\phi=i_{p}^{\ast}f$.
Furthermore, Pukanszky's condition is equivalently satisfied by $\frak
a$. Applying Proposition {\coadjpukan}, the following result is
immediate:

\math{Corollary.}{\semprodcoadjpukan}{\sl A necessary and sufficient
condition for the symplectic subbundle $\eurm F=\big(G\times(\en
+\frak e^{\circ})\big)/D\subset T^{\sharp}(G/D)$ to be
symplectomorphic to the coadjoint orbit ${\script O}^{G}_{\en}$,
is that the symplectic subbundle $\eurm
F_{p}=\big(K_{p}\times(\phi+i_{p}^{\ast}\frak q^{\circ})\big)/P
\subset T^{\sharp}(K_{p}/P)$ be symplectomorphic to the
little-group coadjoint orbit ${\script O}_{\phi}^{K_{p}}$.}

\vskip0.3cm

We determine now an equivalence class in the quotient $\eurm F
=\big(G\times(\en+\frak e^{\circ})\big)/D$, given that $\frak
e^{\circ}=(\frak
q\opluss{\rho}V)^{\circ}\cong\frak
q^{\circ}\subset\k^{\ast}\times\{0\}$ and that $D$ is described by the
exact sequence
$$0\longrightarrow V\buildrel i\over\longrightarrow D\buildrel
j\over\longrightarrow P\longrightarrow e.\eqn\Dexactseq$$
Let $(g,\en+w)\an G\times(\en+\frak e^{\circ})$ and 
$[(g,\en+w)]$ be its equivalence class. We know that $\eurm F$ is an 
affine bundle associated to the principal fibre bundle 
$G\rightarrow G/D$. So, it is sufficient to find the
equivalence class $[g]$. But thanks to the exact sequence
$\Dexactseq$, we obtain $[g]=g\cdot D=[\ek]=\ek\cdot P$ 
if $g=(\ek,v)\an G$ and we may write $\eurm F=
\big(K\times(f+\frak q^{\circ})\big)/P$. Furthermore, there is a
canonical inclusion $K_{p}/P\hookrightarrow K/P$ 
induced by the inclusion of the closed subgroup $K_{p}$ in
$K$; we have so a projection $T_{x}^{\ast}(K/P)
\rightarrow T_{x}^{\ast}(K_{p}/P)$, for $x=[\ek]$,
$\ek\an K_{p}$. 

Let us examine in more detail this projection of
cotangent spaces. If we denote $T_{e}L_{\ek}(\frak p)\subset T_{\ek}K$
as $\frak p_{\ek}$ (recall that $\frak p=\frak a\cap\k$), 
then clearly $T_{x}(K/P)\cong T_{\ek}K/\frak
p_{\ek}$ and $T^{\ast}_{x}(K/P)\cong\frak
p_{\ek}^{\circ}\subset\k_{\ek}^{\ast}$, the
annihilator of $\frak p_{\ek}$. Similarly, $T_{x}(K_{p}/P)\cong
T_{\ek}K_{p}/\frak p_{\ek}$ and $T^{\ast}_{x}(K_{p}/P)$ is isomorphic
to the space of elements of $T^{\ast}K_{p}$ which vanish on
$\frak p_{\ek}$. 
If $T_{\ek}i_{p}\Colon T_{\ek}K_{p}\hookrightarrow
T_{\ek}K$ is the natural inclusion, we may write
$T^{\ast}_{\ek}(K_{p}/P)\cong(T_{\ek}i_{p})^{\ast}\frak p_{\ek}^{\circ}$.
Now, the inclusion $\frak p\subset\frak q$ implies $\frak
q_{\ek}^{\circ}\subset\frak p_{\ek}^{\circ}$; but $\frak
q_{\ek}^{\circ}$ and $(T_{\ek}i_{p})^{\ast}\frak q_{\ek}^{\circ}$ are
the fibres of $\eurm F\cong{\script O}^{G}_{\en}$ and $\eurm
F_{p}\cong{\script O}^{K_{p}}_{\phi}$, respectively, over $x$ (see
also Corollary {\semprodcoadjpukan}). Thus, the fibres of $\eurm
F_{p}$ are obtained from those of $\eurm F$ under the projections
$(T_{\ek}i_{p})^{\ast}$. Therefore:

\math{Corollary.}{\semprodcoadjpukansuppl}{\sl Validity of
Pukanszky's condition for a polarization $\frak h=\frak
a\opluss{\rho}V^{\complex}$ at $\en\an\frak g^{\ast}$ implies that the
coadjoint orbit ${\script O}^{G}_{\en}$ is symplectomorphic to the
quotient $\eurm F=\big(K\times(f+\frak q^{\circ})\big)/P$ and
that the coadjoint orbit 
${\script O}^{K_{p}}_{\phi}$ is obtained by restricting 
$\eurm F\cong{\script O}^{G}_{\en}$ to the closed subset
$K_{p}/P\subset K/P$ and projecting its fibres by the
natural projection between the corresponding cotangent bundles.}

\vskip0.3cm

Clearly, under the conditions of Proposition {\coadjpukansuppl}, the
coadjoint orbit ${\script O}_{\phi}^{K_{p}}$ has properties
analogous to those described in this proposition because, according to
Theorem {\semprodpuk}, Pukanszky's condition on $\frak h$ is
equivalent to the same condition on the polarization $\frak a$.

\chapter{Connections and symplectic induction by semidirect products}

We consider in this section the problem of the arbitrary choice of the
connection in the symplectic induction process, pointed out in {\det}.
This connection plays the role of a Yang-Mills potential in the
more general geometrical interaction scheme due to Guillemin-Sternberg 
{\guist} {\stern} and Weinstein {\wein}.
The symplectic induction is a special case of this model and the
arbitrariness or not of the connection has been related to the 
localization of relativistic particles in references {\de} and {\det}.

We restrict our attention to the case where we have a principal fibre
bundle $\epi\Colon G\rightarrow G/D$, formed by a
semidirect product $G=K\timess{\rho}V$ and a closed
subgroup $D$ belonging to the set of
extensions of the Lie subgroup $P\subset K$ 
by the vector group $V$. In this case, the base space $G/D$ 
is equal to $K/P$; we denote this quotient by $\mit\Sigma$. The following
commutative diagram illustrates this situation.

\vskip0.5cm
\centerline{\beginpicture
\setcoordinatesystem units <1.00000cm,1.00000cm>
\setplotsymbol ({\fiverm .})

\circulararc 180.000 degrees from 10.160 18.256 center at 10.09 18.256
\circulararc 180.000 degrees from 12.700 18.256 center at 12.63 18.256
\circulararc 180.000 degrees from 12.21 16.89 center at 12.21 16.96

\setlinear

\arrow <1mm> [0.1,1] from 5.397 18.7 to  7.144 18.7
\arrow <1mm> [0.1,1] from 7.938 18.7 to  9.684 18.7
\arrow <1mm> [0.1,1] from 10.478 18.7 to 12.224 18.7
\arrow <1mm> [0.1,1] from 13.018 18.7 to 14.764 18.7
\arrow <1mm> [0.1,1] from 5.397 16.8 to  7.144 16.8
\arrow <1mm> [0.1,1] from 7.938 16.8 to  9.684 16.8
\arrow <1mm> [0.1,1] from 13.018 16.8 to 14.764 16.8
\arrow <1mm> [0.1,1] from 10.160 16.510 to 10.160 15.240
\arrow <1mm> [0.1,1] from 12.700 16.510 to 12.700 15.240
\arrow <1mm> [0.1,1] from 10.160 18.256 to 10.160 17.145
\arrow <1mm> [0.1,1] from 12.700 18.256 to 12.700 17.145
\arrow <1mm> [0.1,1] from 10.478 16.701 to 12.224 16.701
\arrow <1mm> [0.1,1] from 12.192 16.891 to 10.446 16.891
\putrule from 10.795 14.926 to 12.065 14.926
\putrule from 10.795 14.806 to 12.065 14.806
\putrule from 7.534  18.400 to 7.534 17.070
\putrule from 7.660  18.400 to 7.660 17.070

\put{$V$} [lB] at  7.458 16.669
\put{$V$} [lB] at  7.458 18.574
\put{$e$} [lB] at 14.9 18.61
\put{$e$} [lB] at 14.9 16.71
\put{$0$} [lB] at  5 16.7
\put{$0$} [lB] at  5 18.6
\put{$G$} [lB] at 10.001 16.669
\put{$D$} [lB] at 10.001 18.574
\put{$K$} [lB] at 12.541 16.669
\put{$P$} [lB] at 12.541 18.574
\put{$G/D$} [lB] at  9.842 14.764
\put{$K/P$} [lB] at 12.383 14.764
\put{$i$} [lB] at 11.303 17.05
\put{$\epi_{1}$} [lB] at 11.2 16.35
\put{$\epi$} [lB] at  9.8 15.843
\put{$\epi_{\frak p}$} [lB] at 12.85 15.875
\put{$=\mit\Sigma$} [lB] at 13.24 14.764
\endpicture
}
\vskip-3.3cm

Let then $\ec\an\Omega^{1}(G)\otimes\don$ be a connection on
$\epi\Colon G\rightarrow G/D$, where $\don$, the Lie
algebra of $D$, is an extension of the Lie algebra $\frak p$ by
$V$. Set now $g=(\ek,v)\an G$, $h=(\el,u)\an D$ and $\xi=(B,b)\an\don$.
Then, if $R_{h}(g)=gh$ is the right action of $D$ on $G$, by the defining
properties of a connection form we must have 
$R_{h}^{\ast}\ec=\Ad(h^{-1})\comp\ec$ and
$\ec_{g}\big(\wh{\xi}(g)\big)=\xi$, where
$\wh{\xi}$ is the fundamental vector field of the Lie
algebra element $\xi$ for the right action on $G$: 
$\wh{\xi}(g)=(\wh{B}(\ek),\ek\cdot b)$. Using the fact that the
$\frak p$-components of $TR_{h}$ and $\Ad(h^{-1})$ are $TR_{\el}$ and
$\Ad(\el^{-1})$ respectively, we obtain immediately that
the $\frak p$-component $\ea=(i^{\ast}\otimes T_{e}\epi_{1})\comp\ec$ 
of the pull-back of $\ec$ under the inclusion $i\Colon
K\hookrightarrow G$, is a connection 1-form on $K\rightarrow K/P$.

Conversely, suppose that we have a connection
$\ea\an\Omega^{1}(K)\otimes\frak p$ on $K\rightarrow K/P$, 
then, for each $\ek\an K$, we
have a horizontal subspace $H_{\ek}\subset T_{\ek}K$ setting
$H_{\ek}=\ker\ea_{\ek}$; $H_{\ek}$ is isomorphic to $T_{q}{\mit\Sigma}$ 
under
the projection $\epi_{\frak p}\Colon K\rightarrow\mit\Sigma$, 
$q=\epi_{\frak
p}(\ek)=[\ek]$. If now $g=(\ek,v)\an G$, one can define a subspace
$\bar{H}_{g}\subset T_{g}G\cong T_{\ek}K\oplus V$, complementary to the
vertical subspace at $g$ and isomorphic also to $T_{q}\mit\Sigma$,
$q=\epi(g)=[g]=[\ek]$. In fact, if $X\an H_{\ek}$, let us define 
$\bar{X}\an\bar{H}_{g}$ by
$$\bar{X}=\big(X,T_{e}\rho(T_{\ek}R_{\ek^{-1}}(X))v\big).\eqn\horvect$$ 
We set for
simplicity $T_{e}\rho(T_{\ek}R^{-1}_{\ek}(X))v=R_{\ek^{-1}}X\cdot v$.
In order to prove that $\bar{H}$ defines a connection on $\epi\Colon
G\rightarrow\mit\Sigma$, 
it is sufficient to check if $\bar{H}$ is invariant
under the right action of $D$ on $G$. Taking $h=(\el,u)\an D$, we
find: $TR_{h}(\bar{X})=(\wt{X},w)$, where $\wt{X}=TR_{\el}(X)$ and
$w=T_{\ek}\rho(X)u+R_{\ek^{-1}}X\cdot v$. But easy calculation shows
that $w=R_{(\ek\el)^{-1}}\wt{X}\cdot(\ek\cdot u+v)$; consequently,
$TR_{h}(\bar{X})\an\bar{H}_{gh}$, which proves that $\bar{H}$ is
indeed a connection. 

In order to calculate now the corresponding connection form $\ec$, we
use the decomposition $(Z,w)=(X+\wh{B}(\ek),R_{\ek^{-1}}X\cdot v+\ek\cdot
b)$ of an arbitrary tangent vector at $(\ek,v)$ into horizontal and
vertical parts, $(B,b)\an\don$, $X\an\ker\ea_{\ek}$. 
Then $\ec$ must satisfy $\ec(Z,w)=(B,b)$ and if we decompose $\ec$ as
$\ec=(\ec_{1},\Delta)$, then $\ec_{1}=\epi_{1}^{\ast}\ea$
and $\Delta\an\Omega^{1}(G)\otimes V$ is given by:
$$\Delta_{(\ek,v)}(Z,w)=\ek^{-1}\cdot(w-R_{\ek^{-1}}X\cdot
v).\eqn\deltaequation$$
We have thus proved:

\math{Proposition.}{\semconnections}{\sl If
$\ec\an\Omega^{1}(G)\otimes\don$ is a connection form on the principal
bundle $G\rightarrow G/D$, then
$\ea=(i^{\ast}\otimes T_{e}\epi_{1})\comp\ec$,
the $\frak p$-component of the
pull-back of $\ec$ under the inclusion $i\Colon K\hookrightarrow G$, is a
connection form on $\epi\Colon K\rightarrow\mit\Sigma$. Furthermore, a
connection $\ea\an\Omega^{1}(K)\otimes\frak p$ determines a preferred 
connection $\ec\an\Omega^{1}(G)\otimes\don$. Indeed,
$\ec$ is equal to $(\epi_{1}^{\ast}\ea,\Delta)$, where $\Delta$ is given
by equation $\deltaequation$ and the horizontal spaces defined by 
$\ec$ are given by equation $\horvect$.}

\math{Remark.}{\indremark}{\it Suppose that $G$ and $D$ 
are such that we can apply symplectic induction from a point
$\en_{0}\an\don^{\ast}$. Then we know {\rm{\det}} that
$\ec_{0}=\en_{0}\comp\ec$ is a 1-form invariant on
the fibres of $\epi$ and therefore,
there exists a 2-form $\beta_{0}\an\Omega^{2}({\mit\Sigma})$
such that $\epi^{\ast}\beta_{0}=d\ec_{0}$. This 2-form gives the
modification term of the canonical symplectic structure of the
cotangent bundle $T^{\ast}(G/D)=T^{\ast}\mit\Sigma$. Now if there exists a
canonical connection $\ea\an\Omega^{1}(K)\otimes\frak p$, then we may
choose $\ec$ in a natural way according to Proposition
{\semconnections} and so $\beta_{0}$ is also canonical.}

\vskip0.3cm

We examine now a special case where the choice of the connection  
$\ec$ is guided by supplementary geometrical structures (so we have 
a canonical connection). The following proposition explains then
a result of {\de} concerning the localization procedure of
relativistic particles with non-zero mass. This is precisely the case
of massive Poincar\'e coadjoint orbits (the hyperboloid $m=const.>0$
is a symmetric space unlike the massless case where the light cone is
not).

\math{Proposition.}{\canonicalconnection}{\sl Let $G=
K\timess{\rho}V$ be a semidirect product,
$\frak g=\k\opluss{\rho} V$ its Lie algebra and
$\en=(f,p)\an\frak g^{\ast}$. Suppose there exists a polarization
$\frak a$ of the complexified Lie algebra 
$\k_{p}^{\complex}$ satisfying Pukanszky's condition
and $K_{p}/P$ and ${Z}=K/K_{p}$ 
are symmetric spaces, where $P$ is 
the closed subgroup of $K_{p}$ determined by $\frak a$ (Lemma
{\polarisationlemma}). Then the coadjoint orbit 
${\script O}^{G}_{\en}$ is symplectomorphic to a symplectic
subbundle of a modified cotangent bundle (Proposition {\coadjpukan})
whose symplectic structure is canonical.}

\undertext{\it Proof}. The only thing we have to prove is that the
modified symplectic structure of the cotangent bundle
$T^{\sharp}(G/D)$ is canonical or, equivalently, that the
connection form $\ea$ is canonical.
Since $K_{p}/P$ and ${Z}=K/K_{p}$ 
are symmetric spaces, there exist involutive
automorphisms $I_{p}\Colon K_{p}\rightarrow K_{p}$ and 
$I\Colon K\rightarrow K$ defining the
canonical symmetric space decompositions $\k_{p}=\frak p\oplus\frak m$
and $\k=\k_{p}\oplus\frak n$, where $\frak m$ and $\frak n$ are the
subspaces of $\k_{p}$ and $\k$ respectively corresponding to the
eigenvalue $-1$ of $T_{e}I_{p}$ and $T_{e}I$. 

Using the known property $\Ad(P)\frak m\subset\frak m$ and
$\Ad(K_{p})\frak n\subset\frak n$ of these subspaces,
we obtain a canonical decomposition 
$\k=\k_{p}\oplus\frak n=\frak p\oplus(\frak m\oplus\frak n)$ with
the same property: $\Ad(P)(\frak m\oplus\frak n)\subset\frak
m\oplus\frak n$. By the invariant connection theory {\kob},
there exists on the
principal fibre bundle $\epi_{\frak p}\Colon K\rightarrow\mit\Sigma$ 
a canonical $K$-invariant connection because we have
a subspace $\frak m\oplus\frak n\subset\k$ invariant under the adjoint
action of $P$ and such that $\k=\frak p\oplus(\frak m\oplus\frak
n)$. Then, Proposition {\semconnections} and Remark {\indremark} finish
the proof.\QED

Consider finally the special case $K_{p}\subset K_{f}$ ($\Rightarrow
Y={Z}$). Then,
Propositions {\coadjpukan} and {\semconnections} can be used in order
to consider the result of Proposition {\modifiedcot} from another
point of view. Indeed, in that case we have
$[\,\k,\k_{p}\,]\subset\ker f$, in particular, $[f]\an{\rm
H}^{1}(\k_{p},\R)$. We can thus apply Corollary {\realpol} and
conclude that $\frak h=\k_{p}^{\complex}\opluss{\rho}V^{\complex}$ is
a real polarization of $\frak g^{\complex}$, satisfying Pukanszky's
condition. Now, $P=K_{p}$ 
and ${\mit\Sigma}=Z$. But for the case of real polarizations, the
content of Proposition {\coadjpukan} is essentially that the coadjoint
orbit in question is isomorphic to a modified cotangent bundle
$(T^{\ast}{\mit\Sigma}, d\theta_{\mit\Sigma}
+\tau^{\ast}\beta_{0})$ (see Remark 3.10
of {\det}). Now, according to Proposition {\semconnections} and Remark
{\indremark}, the choice of a connection $\ea\an\Omega^{1}(K)
\otimes\k_{p}$ on $K\rightarrow K/P$ determines this 2-form
completely. The connection
$\ea$ in turn, is determined if we fix a subspace $\frak n\subset\k$
such that $\k_{p}\oplus\frak n=\k$ and $\Ad(K_{p})\frak n=\frak
n$. Then, $\ea$ is the $\k_{p}$-component of the Maurer-Cartan form 
on $K$. It must be emphasized here that, the 2-form
$\alpha_{0}\an\Omega^{2}({Z})$ appearing in Proposition
{\modifiedcot} and giving the modification term of the
symplectic structure of $T^{\ast}{Z}$, depends on the point
$(f,p)$ of the coadjoint orbit for which $K_{p}\subset K_{f}$. On the
other hand, we have just seen that the 2-form $\beta_{0}\an\Omega^{2}
({\mit\Sigma})$
depends on the choice
of a connection on the principal bundle $G\rightarrow G/D$. Thus, the
differential forms $\alpha_{0}$ and $\beta_{0}$ are not canonical and in 
general $\alpha_{0}\neq\beta_{0}$. But in any case,
Proposition {\coadjpukan} tells us that the symplectic structures
$\omega_{Z}+\tau^{\ast}\alpha_{0}$ and $\omega_{\mit\Sigma}
+\tau^{\ast}\beta_{0}$ are equivalent, that is,
there exists a bijection $T^{\ast}{Z}
\rightarrow T^{\ast}\mit\Sigma$ which is a symplectomorphism with
respect to these structures.

\chapter{Symplectic induction and semidirect product}

We have seen previously (Proposition {\coadjpukan}) that the validity
of Pukanszky's condition for a polarization of the coadjoint orbit 
${\script O}^{G}_{\en}$ is equivalent to the fact that this
orbit is symplectomorphic to a subbundle of a 
modified cotangent bundle, obtained by
symplectic induction from a point. In this section we will discuss a
more general property of the coadjoint orbits of a
semidirect product. See {\land} for an equivalent approach.

We state now the principal result of this section.

\math{Theorem.}{\semprodind}{\sl The coadjoint orbit ${\script
O}^{G}_{\en}$ through $\en=(f,p)\an\frak g^{\ast}$
of a semidirect product $G=
K\timess{\rho}V$, is always obtained by symplectic induction from the
coadjoint orbit ${\script O}^{G_{p}}_{\en_{p}}$ of $G_{p}$
passing through $\en_{p}=(i_{p}^{\ast}f,p)\an\frak g_{p}^{\ast}$, 
with groups $G=K\timess{\rho}V$ and $G_{p}=
K_{p}\timess{\rho}V$: $${\script O}_{\en}^{G}=
{\rm Ind}_{G_{p}}^{G}({\script O}^{G_{p}}_{\en_{p}}).$$
Note that ${\script O}^{G_{p}}_{\en_{p}}=G_{p}\cdot\en_{p}$ 
is canonically isomorphic to ${\script
O}^{K_{p}}_{\phi}=K_{p}\cdot\phi$.}

\undertext{\it Proof}. Using the notation of the section on
symplectic induction, let us choose the symplectic manifold $M$ as
$M=G_{p}\cdot\en_{p}$ and the groups $G$
and $H$ as $G=K\timess{\rho}V$ and $H=
K_{p}\timess{\rho}V=G_{p}$. We will then apply symplectic induction
from $M$ with the above mentioned groups.

In our case, one can consider the symplectic manifold $M$ as the
coadjoint orbit of $K_{p}$ passing through
$\phi$, because
$(\ek,v)\cdot(\phi,p)=\big(\ek\cdot\phi+i_{p}^{\ast}(p\odot v),p\big)
=(\ek\cdot\phi,p)$,
$\forall(\ek,v)\an G_{p}$, since $i_{p}^{\ast}(p\odot v)=0$ (see Lemma
{\rangeoftau}). 
The action of $G_{p}$ on $M$ is Hamiltonian
with momentum mapping $J_{M}\Colon M\rightarrow\frak g_{p}^{\ast}$
given by $J_{M}(m)=m$, $m=(\ek\cdot\phi,p)$.

Using the conventions of section 7, one readily verifies that the zero
level set of the momentum map $J_{\check{M}}\Colon\check{M}=M\times
T^{\ast}G\rightarrow\frak h=\frak g_{p}$ is given by
$$J_{\check{M}}^{-1}(0)=\{\big((\varphi,p),g,(z,w)\big)\an M
\times T^{\ast}G\,|\,\varphi=i_{p}^{\ast}z,\,w=p,g\an G\},$$
and knowing that
$\varphi=\ek\cdot\phi=\ek\cdot i_{p}^{\ast}f=i_{p}^{\ast}(\ek\cdot
f)$ for some $\ek\an K_{p}$, we have the following characterization
for $J_{\check{M}}^{-1}(0)$:
$$J_{\check{M}}^{-1}(0)
=\big\{\big((\ek\cdot\phi,p),g,(\ek\cdot f+p\odot
v, p)\big)\,|\,\ek\an K_{p},g\an G,v\an V\big\}$$
Then, direct calculation shows that the point
$\big((\ek\cdot\phi,p),g,(\ek\cdot f+p\odot v,p)\big)$ of
$J_{\check{M}}^{-1}(0)$ can be
represented in the quotient
$J_{\check{M}}^{-1}(0)/G_{p}$  by the point $\big(\hat{\ek}\cdot
(\ek\cdot f+p\odot v)+q\odot\hat{v},q\big)$ if $g=(\hat{\ek},\hat{v})$,
where $q$ represents $g$ in $G/G_{p}=Z$:
$q=\hat{\ek}\cdot p$. 
We realize thus easily that the points of the 
induced manifold $M_{ind}={\rm
Ind}_{G_{p}}^{G}(G_{p}\cdot\en_{p})=J_{\check{M}}^{-1}(0)/G_{p}$ 
will be of the form $\big(\el\cdot f+\el\cdot p\odot u,\el\cdot p\big)$, 
$\el\an K$, $u\an V$,
so $M_{ind}$ and ${\script O}_{\en}^{G}$ are isomorphic as manifolds. 

In order to establish a symplectomorphism between $M_{ind}$ and
${\script O}^{G}_{\en}$, we proceed as follows. The left action of
$G$ on $\check{M}$, obtained by taking the cotangent lift of the left
action of $G$ on itself and letting $G$ act trivially on $M$, projects
on $M_{ind}\cong{\script O}^{G}_{\en}$ as the coadjoint action of $G$,
as one verifies by easy calculation. We choose now an element
$g_{0}=(\ek_{0},v_{0})\an G$ and let
$n_{0}=\big((\phi,p),(\ek_{0},v_{0}),(f,p)\big)\an
J^{-1}_{\check{M}}(0)$. The image of $n_{0}$
in $M_{ind}$ is equal to $(\ek_{0}\cdot f+\ek_{0}\cdot p\odot
v_{0},\ek_{0}\cdot p)=(h,q)$. Furthermore, the projection of a vector
at $n_{0}$ induced by the action of $G$ on $\check{M}$, 
will coincide with the corresponding vector at $(h,q)\an M_{ind}$
induced by the coadjoint action of $G$ on $M_{ind}$. So, consider an
element $\xi=(A,a)\an\frak g$; then $\xi_{\check{M}}(n_{0})=
\big((0,0),T_{e}R_{g_{0}}(\xi),(0,0)\big)$ and  $[\xi_{\check{M}}
(n_{0})]=\xi_{\frak g^{\ast}}(h,q)$. Similarly, if
$\eta=(B,b)\an\frak g$, then:
$$\eqalign{(\omega_{ind})_{(h,q)}(\xi_{\frak g^{\ast}},
\eta_{\frak g^{\ast}}\big)&=-(f,p)\big(\big[
TL_{g_{0}^{-1}}\comp TR_{g_{0}}(A,a),
TL_{g_{0}^{-1}}\comp TR_{g_{0}}(B,b)\big]\big)\cr
\hfil&=-(h,q)([(A,a),(B,b)])\cr}.$$
This shows that $\omega_{ind}$ coincides with the standard symplectic
structure of the coadjoint orbit ${\script O}^{G}_{\en}$. \QED

\vskip0.3cm

We observe here the following analogy with the construction of Rawnsley
{\ra}, described in section 3. According to Lemma {\coadorbitstruc},
we have the fibration 
${\script O}^{G}_{\en}\rightarrow Y
\rightarrow Z$, where the typical fibres are respectively
$p\odot V$ and ${\script O}^{K_{p}}_{\phi}$.
But now Theorems {\indmanifold} and {\semprodind}
ensure that we have indeed a (non canonical)
fibration ${\script O}^{G}_{\en}
\rightarrow T^{\ast}Z$, with typical fibre ${\script
O}^{K_{p}}_{\phi}$.

\math{Remark.}{\trivial}{\it It is evident that if ${Z}$ is a
contractible space, then the coadjoint orbit ${\script O}^{G}_{\en}$
is globally diffeomorphic to the product $T^{\ast}{Z}
\times{\script O}^{K_{p}}_{\phi}$. In particular, when $K_{p}$ is a
group-deformation retract of $K$, that is when there exists a homotopy 
$H\Colon[0,1]\times K
\rightarrow K$ with the properties $H(0,\ek)=i_{p}(r(\ek))$,
$H(1,\ek)=\ek$ and $H_{t}(\ek\el)=H_{t}(\ek)H_{t}(\el)$,
$H_{t}(\ek)=H(t,\ek)$, then $Z$ is contractible; a homotopy 
$\check{H}\Colon[0,1]\times {Z}\rightarrow{Z}$ between
the constant map and the identity on $Z$, is given by
by $\check{H}(t,[\ek])=[H(t,\ek)]$.}

\vskip0.3cm

As an immediate application of Theorem {\semprodind}, we discuss in
the light of symplectic induction the result of Proposition
{\modifiedcot}. If $K_{p}\subset K_{f}$, then ${\script
O}^{K_{p}}_{\phi}=\{\phi\}$ and the coadjoint orbit of $\en$ is hence
obtained by symplectic induction from a point. Thus, according to
Proposition 2.11 of {\det}, ${\script O}^{G}_{\en}$ must be isomorphic
to a modified cotangent bundle
$T^{\sharp}(G/G_{p})=T^{\sharp}Z$, where the modification term
is determined by a connection on the principal bundle $G\rightarrow
G/G_{p}$, in accordance with Proposition {\modifiedcot}. 
The choice of this connection has been discussed in section 9.

\chapter{Examples}

We consider here three representative examples of semidirect product
and we apply the general formalism developed in the previous sections.
The semidirect product Lie groups we analyze below, are important for
the non-relativistic particle dynamics.

{\bf $\bullet$
The special Euclidean group of $\R^{3}$.} Let $K=SO(3)$ be the
Lie group of rotations in $V=\R^{3}$ preserving the standard scalar
product $\langle-,-\rangle\Colon\R^{3}\times\R^{3}\rightarrow\R$. The
familiar representation of the elements of $SO(3)$ as $3\times 3$
matrices, enables us to form the semidirect product
$G=SE(3)=SO(3)\semidir\R^{3}$, the Euclidean group in $\R^{3}$.

The Lie algebra $\frak s\frak e(3)$ as well as its dual 
$\frak s\frak e(3)^{\ast}$ are canonically isomorphic to
$\R^{3}\oplus\R^{3}$. Easy calculation shows that if $p\an
V^{\ast}\cong\R^{3}$, then the linear map $\tau_{p}^{\ast}\Colon
V=\R^{3}\rightarrow\k^{\ast}\cong\R^{3}$ is given by
$\tau_{p}^{\ast}(v)=p\times v$ (the usual cross product of the
vector space $\R^{3}$). 

Now let $\en=(f,p)\an\frak s\frak e(3)^{\ast}$ be an element such that
$f=s\bit u$, $p=k\bit u$ and $\langle \bit u,\bit u\rangle=1$ 
($s,k>0$). Then $K_{p}=K_{f}=(K_{p})_{\phi}\cong SO(2)$, the orbits
$Z$ and $Q$ are 2-spheres $\Sphere^{2}$
and ${\script O}^{K_{p}}_{\phi}=\{\phi\}$. We can furthermore apply 
Propositions {\submanifolds} and {\modifiedcot} which show that the
coadjoint orbit of $\en$ coincides as a manifold to $T^{\ast}\Sphere^{2}$,
but its symplectic structure is modified by a ``spin term" which,
in this case, is $s$-times the canonical symplectic structure of
$\Sphere^{2}$ (its volume element). 

Taking into account the discussion after Proposition
{\canonicalconnection}, we can reconsider this result in the context
of algebraic polarizations: since $[f]\an{\rm H}^{1}(\k_{p},\R)$ ($\k_{p}=
\frak s\frak o(2)$ is an abelian Lie algebra), Corollary {\realpol} can 
be applied; therefore, the subspace $\frak
a=\frak s\frak o(2)^{\complex}$ is a real polarization of
$\k_{p}^{\complex}=\frak s\frak o(2)^{\complex}$ satisfying
Pukanszky's condition, so $$\frak h=\frak s\frak
o(2)^{\complex}\oplus(\R^{3})^{\complex}$$
is a real polarization of $\frak g^{\complex}$ (with respect to $\en$)
satisfying also the same condition (see also Proposition {\semprodpol}
and Theorem {\semprodpuk}). Then, Proposition {\coadjpukan} gives the
same result on the structure of ${\script O}^{G}_{\en}$.
Alternatively, one could use Theorem {\semprodind} (see discussion
after Remark {\trivial}). 

For this polarization, the
groups $D$ and $E$ coincide with $G_{p}=SO(2)\semidir\R^{3}$.
Furthermore, there exists a canonical connection on the principal bundle
$G\rightarrow G_{p}$ because the spaces $K_{p}/P=\hbox{point}$ and
$K/K_{p}=\Sphere^{2}$ are symmetric spaces (see Proposition
{\canonicalconnection}).

{\bf $\bullet$
The Galilei group of $\R^{3}\oplus\R$.} Take now as group $K$ the
Lie group $SE(3)$ of the previous example, $K=SO(3)\semidir\R^{3}$ and
the vector space $V$ as $\R^{3}\oplus\R$ (Galilean space-time). We
have a representation $\rho\Colon K\rightarrow GL(V)$ given by
$$\rho(R,\bit b)=\left(\matrix{R&\bit b\cr
            		  0&1}\right)$$
and consequently one can consider the semidirect product $G=K\semidir
V$. We recognize $G$ as the Galilei group in dimension 3+1, see
{\sou}.

Using the isomorphism $\frak g\cong\frak
g^{\ast}\cong(\R^{3}\oplus\R^{3})\oplus(\R^{3}\oplus\R)$, we may
represent the elements $f\an\k^{\ast}$ and $p\an V^{\ast}$ as
$f=(\bit l,\bit g)$, $\bit l,\bit g\an\R^{3}$ and $p=(\bit p,E)$,
$\bit p\an\R^{3},E\an\R$. Under these identifications, if
we set $\ek=(R,\bit b)\an SE(3)$ and 
$x=\left(\matrix{\bit r\cr t}\right)\an\R^{3}\times\R$, one readily finds:
$$p\odot x=(\bit p\times \bit r,\bit pt)\quad\hbox{and}\quad\ek\cdot
p=(R\bit p,E-\langle R\bit p,\bit b\rangle).
\eqn\contraodot$$
\indent $(i)$ Let us choose $\en=(f,p)\an\frak g^{\ast}$ as $f=
(s\bit u,0)$, $p=(k\bit u,E)$, $s,k>0$, $\langle \bit u,\bit
u\rangle=1$. This choice corresponds to the standard non-relativistic
particle of mass zero with spin $s$ and color $k$, see {\sou}. 

By formula $\contraodot$ one easily finds
$K_{p}=SO(2)\semidir\R^{2}$, $\R^{2}$ being the subspace perpendicular
to $\bit u$ and $SO(2)$ the rotation group of this
subspace. In this case, only the $\frak s\frak o(2)^{\ast}$-component
of $\phi=i_{p}^{\ast}f\an\k_{p}^{\ast}$ is non-zero and consequently,
${\script O}^{K_{p}}_{\phi}$ is a coadjoint orbit of $SO(2)$; so 
${\script O}^{K_{p}}_{\phi}=\{\phi\}$ because $SO(2)$ is abelian.
Furthermore, the homogeneous space ${Z}=K/K_{p}$ is simply
the product $\Sphere^{2}\times\R$. Thus, according to
Theorem {\semprodind} (see also discussion after Remark {\trivial}), 
the coadjoint orbit of $\en$ is symplectomorphic
to a modified cotangent bundle $T^{\sharp}(K/K_{p})=T^{\sharp}\Sphere^{2}
\times\R^{2}$.

One can obtain the same result using the technique of polarizations.
Indeed, with the previous choices, we have $[f]\an{\rm H}^{1}(\k_{p},\R)$
(because $SO(2)$ is abelian) and so Corollary {\realpol} can be
applied. The real polarization $\frak h$ provided by this corollary is
$$\frak h=\big(\frak s\frak
o(2)^{\complex}\oplus(\R^{2})^{\complex}\big)\oplus(\R^{3}\oplus
\R)^{\complex}$$
and the groups $D$ and $E$
$$D=E=(SO(2)\semidir\R^{2})\semidir(\R^{3}\times\R).$$
Then, according to Proposition {\coadjpukansuppl}, the coadjoint orbit
${\script O}^{G}_{\en}$ is symplectomorphic to a modified cotangent
bundle $T^{\sharp}(K/K_{p})=T^{\sharp}\Sphere^{2}
\times\R^{2}$.

$(ii)$ We choose now an element $\en=(f,p)\an\frak g^{\ast}$ 
setting $p=(\bit p,0)$, $\langle \bit p,\bit p\rangle=1$ 
and $f=(0,\bit g)$ with 
$\langle \bit g,\bit p\rangle=0$ and $\langle \bit g,\bit g\rangle=1$.
The coadjoint orbit of $\en$ has a less evident interpretation;
according to {\guist}, it corresponds to ``particles at infinity with
infinite velocity and mass zero".

Now, by equation $\contraodot$ we find $K_{p}=SO(2)\semidir\R^{2}$,
where $\R^{2}$ is the 2-dimensional subspace of $\R^{3}$ perpendicular
to the line $\R\bit p$ and $SO(2)$ the special orthogonal group of
this subspace. On the other hand, the projection $\phi=i_{p}^{\ast}f$
is equal to $f$ and we readily obtain $(K_{p})_{\phi}=\{e\}\semidir\R
\bit g$. Furthermore, $\ker\tau_{p}^{\ast}=\R\bit p$, as one can see
from $\contraodot$. Thus, the isotropy subgroup $G_{\en}$ is
2-dimensional and the orbit ${\script O}^{G}_{\en}$ will be
8-dimensional. Indeed, as in the previous example,
the orbit ${Z}=K/K_{p}$ is the product $\Sphere^{2}\times\R$
and by Theorem 
{\semprodind}, the coadjoint orbit of $\en$ can be identified (in a
non canonical way) with a fibre bundle over $T^{\ast}{Z}=
T^{\ast}\Sphere^{2}\times\R^{2}$ whose
typical fibre is ${\script O}^{K_{p}}_{\phi}=K_{p}/(K_{p})_{\phi}\cong
T^{\ast}\Sphere^{1}\cong\Sphere^{1}\times\R$.

However, we can further investigate the structure of this orbit as
follows. Observe first that the subspace $\frak
a=(0\oplus\R^{2})^{\complex}\subset\k_{p}^{\complex}$ is a real
polarization with respect to $\phi$. In fact,
$(\k_{p})_{\phi}^{\complex}\subset\frak a$, $\frak a$ is invariant
under the adjoint action of $(K_{p})_{\phi}$ (see relation
$\adjoint$), $\dim_{\complex}\frak a={1 \over
2}(\dim\k_{p}+\dim(\k_{p})_{\phi})$ and $[\,\frak a,\frak a\,]=0$. In
this case $\don=\frak e=0\oplus\R^{2}$ and $D=E=\{e\}\semidir\R^{2}$.
As a result, $D\cdot\phi=(\R,\bit g)=\phi+\frak e^{\circ}$, which
means that $\frak a$ satisfies Pukanszky's condition. By Theorem 
{\semprodpuk}, the subspace
$$\frak h=\frak a\oplus V^{\complex}=(0\oplus\R^{2})^{\complex}\oplus
(\R^{3}\oplus\R)^{\complex}\subset\frak g^{\complex}$$
is a real polarization of $\frak g^{\complex}$ (with respect to $\en$)
satisfying also Pukanszky's condition. Then Proposition
{\coadjpukansuppl},
applied for a real polarization, tells us that the coadjoint orbit
${\script O}^{G}_{\en}$ is symplectomorphic to a modified cotangent
bundle $T^{\sharp}(G/D)\cong T^{\sharp}(SO(3)\times\R)$. In particular,
${\script O}^{G}_{\en}\cong SO(3)\times\R^{5}$ as a manifold.

{\bf $\bullet$
The Bargmann group of $\R^{3}\oplus\R$.} Consider again the special
euclidean group $SE(3)$ of $\R^{3}$ and let $\rho\Colon
SE(3)\rightarrow GL(\R^{5})$ be the representation given by
$$\rho(R,b)=\left(\matrix{R&\bit b&0\cr
                          0&1&0\cr
                          -\bit b^{t}R&-\bit b^{2}/2&1}\right).$$
The semidirect product $G=SE(3)\timess{\rho}\R^{5}$ is called
Bargmann group and it is a non-trivial extension of the Galilei
group, previously studied, by $\R$, see {\dbkp}.
If we write an element $p\an\R^{5}$
as $p=(\bit p,E,m)$ and $\ek=(R,\bit b)\an K=SE(3)$, we find easily
$$\ek\cdot p=(R\bit p+m\bit b,E-\langle R\bit p,\bit b\rangle-m{
\bit b^{2}\over 2},m).\eqn\bargcontr$$
Let us consider the coadjoint orbit of the element $\en=(f,p)\an\frak
g^{\ast}$ with $f=(s\bit u,0)$ and $p=(0,0,m)$, $m>0$, $s$ and $\bit
u$ as above. Easy calculation using $\bargcontr$ shows that
$K_{p}=SO(3)\semidir\{0\}$ and consequently the projection
$i_{p}^{\ast}\Colon\k^{\ast}\rightarrow\k_{p}^{\ast}$ is simply the
projection on to the first factor, $i_{p}^{\ast}f=\phi=s\bit u$. Now the
orbit ${Z}=K/K_{p}$ is simply ${Z}=\R^{3}$ and
${\script O}^{K_{p}}_{\phi}\cong\Sphere^{2}$, thus by Remark {\trivial}
the coadjoint orbit ${\script O}^{G}_{\en}$ of $\en$ is diffeomorphic
to the product $T^{\ast}\R^{3}\times\Sphere^{2}$. We recognize here
the phase space of non-relativistic particles of mass $m$ and spin
$s$, {\sou}.

One could also investigate the structure of ${\script O}^{G}_{\en}$
using the technique of algebraic polarizations. To this end, one
proceeds as follows. Consider the subspaces
$$\frak a^{\pm}_{0}=\{A\an\k_{p}^{\complex}\;|\;\hat{A}\times u=\pm
i\hat{A}\}$$
of $\k_{p}^{\complex}=\frak s\frak o(3)^{\complex}$, where
$A\mapsto\hat{A}$ is the natural isomorphism $\frak s\frak
o(3)^{\complex}\cong(\R^{3})^{\complex}$. These subspaces have
complex dimension one and are such that $\frak a^{+}_{0}\oplus\frak
a^{-}_{0}\oplus(\k_{p})_{\phi}^{\complex}=\k_{p}^{\complex}$. Furthermore,
it is elementary to verify that $[\,\frak a^{+}_{0},\frak
a^{+}_{0}\,]\subset(\k_{p})_{\phi}^{\complex}$ and $[\,\frak
a^{+}_{0},(\k_{p})_{\phi}^{\complex}]\subset\frak a^{+}_{0}$ 
(similarly for $\frak a^{-}_{0}$). Thus, if we set 
$$\frak a^{\pm}
=\frak a^{\pm}_{0}\oplus(\k_{p})_{\phi}^{\complex},$$
we obtain two
(isomorphic) complex subalgebras of $\k_{p}^{\complex}$ such that complex
conjugation on the one yields the other. This means that $\frak
a^{\pm}$ are (isomorphic) complex polarizations of $\k_{p}^{\complex}$
with respect to $\phi=s\bit u$. The real Lie subalgebras $\don$ and $\frak
e$ of $\k_{p}$ are easily found to be
$\don=(\k_{p})_{\phi}=\frak s\frak o(2)$ and $\frak e=\k_{p}=\frak
s\frak o(3)$. It is then evident that $D\cdot\phi=\phi$ and so $\frak
a^{\pm}$ satisfy Pukanszky's condition because $\frak e^{\circ}=0$. 
As a result, the subalgebras
$$\frak h^{\pm}=\frak a^{\pm}\opluss{\rho}(\R^{5})^{\complex}$$
of $\frak g^{\complex}$ are complex polarizations satisfying
Pukanszky's condition. Therefore, by Proposition 
{\coadjpukansuppl} we conclude that the coadjoint orbit
${\script O}^{G}_{\en}$ is a fibre bundle over
$T^{\ast}(G/E)=T^{\ast}\R^{3}$ with typical fibre $E/D=\Sphere^{2}$.

\vskip1cm
   \ifreferenceopen \Closeout\referencewrite \referenceopenfalse \fi
   \line{\bf\hskip0pt\hfil References\hfil}\vskip\headskip
   \vskip0.3cm
   \catcode `\^^M=10
\refitem {\tenprm [\tenpbf 1\tenprm ]}

\author {Duval, C., Burdet, G., K\"unzle, H. P., Perrin, M.}
\titlosa {Bargmann structures and Newton-Cartan theory}
\periodiko {Phys. Rev. D}
\volume {31(8)}
\selides {1841--1853 (1985)}
\refitem {\tenprm [\tenpbf 2\tenprm ]}

\author {Duval, C., Elhadad, J.}
\titlosa {Geometric quantization and localization of relativistic spin^^Msystems}
\periodiko {Contemp. Math.}
\volume {132}
\selides {317--330 (1992)}
\refitem {\tenprm [\tenpbf 3\tenprm ]}

\author {Duval, C., Elhadad, J., Tuynman, G.M.}
\titlosa {Pukanszky's condition and symplectic induction}
\periodiko {J. Diff. Geom.}
\volume {36}
\selides {331--348 (1992)}
\refitem {\tenprm [\tenpbf 4\tenprm ]}

\author {Guillemin, V., Sternberg, S.}
\titlosb {Symplectic techniques in physics}
\ekdoths {Cambridge University Press, Cambridge (1991)}
\refitem {\tenprm [\tenpbf 5\tenprm ]}

\author {Kazhdan, D., Kostant, B., Sternberg, S.}
\titlosa {Hamiltonian Group Actions and Dynamical Systems of Calogero^^MType}
\periodiko {Comm. Pure Appl. Math.}
\volume {31}
\selides {481--508 (1978)}
\refitem {\tenprm [\tenpbf 6\tenprm ]}

\author {Kobayashi, S., Nomizu, K.}
\titlosb {Foundations of Differential Geometry}
\ekdoths {Vol. I, II, Interscience Publishers,^^MJohn Wiley \& Sons, (1963, 1969)}
\refitem {\tenprm [\tenpbf 7\tenprm ]}

\author {Kostant, B.}
\titlosa {On certain unitary representations which arise from a^^Mquantization theory}
\periodiko {in ``Group representations in mathematics and physics"}
\selides {(V. Bargman, ed.), pp. 237--253,^^MProceedings Battelle Seattle 1969 Rencontres; LNP {\tenpbf 6},}
\selides {Springer-Verlag, Berlin, 1970}
\refitem {\tenprm [\tenpbf 8\tenprm ]}

\author {Landsman, N. P.}
\titlosa {Rieffel induction as generalized quantum Marsden-Weinstein^^Mreduction}
\periodiko {J. Geom. Phys.}
\volume {15}
\selides {285--319 (1995)}
\refitem {\tenprm [\tenpbf 9\tenprm ]}

\author {Marsden, J. E., Ratiu, T., Weinstein, A.}
\titlosa {Semidirect products and reduction in mechanics}
\periodiko {Trans. A.M.S}
\volume {281}
\selides {147--177 (1984)}
\refitem {\tenprm [\tenpbf 10\tenprm ]}

\author {Pukanszky, L.}
\titlosa {On the theory of exponential groups}
\periodiko {Trans. A.M.S}
\volume {126}
\selides {487--507 (1967)}
\refitem {\tenprm [\tenpbf 11\tenprm ]}

\author {Rawnsley, J. H.}
\titlosa {Representations of a semidirect product by quantization}
\periodiko {Math. Proc. Camb. Phil. Soc.}
\volume {78}
\selides {345--350 (1975)}
\refitem {\tenprm [\tenpbf 12\tenprm ]}

\author {\'Sniatycki, J.}
\titlosb {Geometric quantization and quantum mechanics}
\ekdoths {Appl. Math. Sci.,}
\volume {50}
\ekdoths {Springer-Verlag (1980)}
\refitem {\tenprm [\tenpbf 13\tenprm ]}

\author {Souriau, J.-M.}
\titlosb {Structure des syst\`emes dynamiques}
\ekdoths {Dunod, Paris (1969)}
\refitem {\tenprm [\tenpbf 14\tenprm ]}

\author {Sternberg, S.}
\titlosa {Minimal coupling and the symplectic mechanics of a classical^^Mparticle in the presence of a Yang-Mills field}
\periodiko {Proc. Nat. Acad. Sci.}
\volume {74}
\selides {5253--5254 (1977)}
\refitem {\tenprm [\tenpbf 15\tenprm ]}

\author {Weinstein, A.}
\titlosa {A universal phase space for particles in Yang-Mills fields}
\periodiko {Lett. Math. Phys.}
\volume {2}
\selides {417--420 (1978)}
\refitem {\tenprm [\tenpbf 16\tenprm ]}

\author {Zakrzewski, S.}
\titlosa {Induced representations and induced Hamiltonian actions}
\periodiko {J. Geom. Phys.}
\volume {3}
\selides {211--219 (1986)}
\refitem {\tenprm [\tenpbf 17\tenprm ]}

\author {Ziegler, F.}
\titlosb {M\'ethode des orbites et repr\'esentations quantiques}
\ekdoths {Th\`ese, Univesrit\'e de Provence, D\'ecembre 1996}
\catcode `\^^M=5


\end